\documentclass[10pt,a4paper,reqno,smallextended]{amsart}
\usepackage{a4wide}
\usepackage{latexsym}
\usepackage{enumitem}
\usepackage{amsmath}
\usepackage{amssymb,textcomp}
\DeclareFontFamily{OT1}{pzc}{}
\DeclareFontShape{OT1}{pzc}{m}{it}%
            {<-> s * [1.100] pzcmi7t}{}
\DeclareMathAlphabet{\mathscr}{OT1}{pzc}{m}{it}
\usepackage{geompsfi-pdf}

\renewcommand{\d}{\mathrm d}

\newcommand{\rme}{\mathrm e}
\newcommand{\rmS}{\mathrm S}
\newcommand{\cA}{\mathcal A}
\newcommand{\A}{\mathcal A}
\newcommand{\cJ}{\mathcal J}
\newcommand{\J}{\mathcal J}
\newcommand{\sfd}{\mathsf d}

\newcommand{\gradG}{\nabla_{\!G}}
\newcommand{\Tan}[2]{\mathrm T_{#1}(#2)}
\newcommand{\CoTan}[2]{\mathrm T_{#1}^*(#2)}
\newcommand{\dm}[2]{\d#1}
\newcommand{\dmm}[2]{#1(\d#2)}
\newcommand{\CE}[3]{\mathsf{CE}(#1,#2;#3)}
\newcommand{\gdens}{g_\e}
\newcommand{\hgdens}{\hat g_\e}
\newcommand{\Md}[2]{M(#1,#2)}

\newtheorem{lemma}{Lemma}[section]
\newtheorem{definition}[lemma]{Definition}

\newtheorem{theorem}[lemma]{Theorem}

\numberwithin{equation}{section}

\newcommand{\R} {{\mathbb R}}

\def\Z{\mathcal Z}

\def\L{\mathscr L}
\def\M{\mathscr M}

\let\e\varepsilon

\let\epsilon\varepsilon
\def\F{\mathcal F}
\def\E{\mathcal E}

\def\hE{\hat \E}
\def\hJ{\hat \J}
\def\hue{\hat u_\e}
\def\huz{\hat u_0}

\def\hge{\hat \gamma_\e}
\def\hgz{\hat \gamma_0}

\let\weakto\rightharpoonup
\def\weakstarto{\stackrel*\rightharpoonup}

\DeclareMathOperator\supp{supp}
\DeclareMathOperator\Tr{Tr}
\def\longrightharpoonup{\relbar\joinrel\rightharpoonup}
\def\longweakstarto{\stackrel*\longrightharpoonup}
\newcommand{\restr}[1]{\hbox{$|_{#1}$}}
\def\mass{m}
\newenvironment{remark}%
  {\par\medbreak\refstepcounter{lemma}%
    \noindent\textbf{Remark~\thelemma. }}%
  {\qed\par\medskip}

\usepackage{bibentry}

\begin{document}

\title{Passing to the Limit in a Wasserstein Gradient Flow: From
 Diffusion to Reaction} 

\author[Arnrich, Mielke, Peletier, Savar\'e,
and Veneroni]{Steffen Arnrich$^1$ \and Alexander Mielke$^2$ \and Mark
 A. Peletier$^3$ \and Giuseppe Savar\'e$^4$ \and Marco Veneroni$^5$}
 
\footnotetext[1]{Department of Mathematics and Computer Sciences, Technische Universiteit Eindhoven}
\footnotetext[2]{Weierstra\ss-Institut f\"ur Angewandte
Analysis und Stochastik and Humboldt Universit\"at zu Berlin}
\footnotetext[3]{Institute for Complex Molecular Systems and Department of Mathematics and Computer Sciences, Technische Universiteit Eindhoven; \url{m.a.peletier@tue.nl}} 
\footnotetext[4]{Dipartimento di Matematica ``F.~Casorati'', Universit\`a di Pavia}
\footnotetext[5]{Department of Mathematics and Statistics, McGill University, Montreal}

\date{\today}

\begin{abstract}
% Part of the bibentry package
\nobibliography*
We study a singular-limit problem arising in the modelling of chemical reactions. At finite $\e>0$, the system is described by a Fokker-Planck convection-diffusion equation with a double-well convection potential. This potential is scaled by $1/\e$, and in the limit $\e\to0$, the solution concentrates onto the two wells, resulting into a limiting system that is a pair of ordinary differential equations for the density at the two wells. 

This convergence has been proved in~Peletier, Savar\'e, and Veneroni,  {\em SIAM Journal on Mathematical Analysis}, 42(4):1805--1825, 2010, using the linear structure of the equation. In this paper we re-prove the result by using solely the Wasserstein gradient-flow structure of the system. In particular we make no use of the linearity, nor of the fact that it is a second-order system. 

The first key step in this approach is a reformulation of the equation as the minimization of an action functional that captures the property of being a  \emph{curve of maximal slope} in an integrated form. The second important step is a rescaling of space. Using only the Wasserstein gradient-flow structure, we prove that the sequence of rescaled solutions is pre-compact in an appropriate topology. We then prove a Gamma-convergence result for the functional in this topology, and we identify the limiting functional and the differential equation that it represents. A consequence of these results is that solutions of the $\e$-problem converge to a solution of the limiting problem. 
\end{abstract}

\maketitle


\tableofcontents

\section{Introduction}
\label{sec:introduction}

In a seminal paper in 1940, Kramers introduced a model of chemical reactions in which the system is represented by a Brownian particle in a potential energy landscape~\cite{Kramers40}. In this model the wells of the potential energy correspond to stable states of the system, and a reaction event is the passage of the particle from one well to another. 
By analyzing the probability of such a reaction event in terms of system parameters, Kramers was able to improve existing formulas for the macroscopically observed reaction rate. 

Although Kramers does not state it in these terms, the central result in~\cite{Kramers40} is a convergence result in the limit of \emph{large activation energy}. 
In~\cite{PeletierSavareVeneroni10} we provided a first rigorous proof of this result in the case of Brownian particles without inertia. The present paper can be considered a sequel to~\cite{PeletierSavareVeneroni10}, in which we address a question that was left unanswered in~\cite{PeletierSavareVeneroni10}. 

The issue hinges on the fact that the system of~\cite{PeletierSavareVeneroni10} is a gradient flow of a free-energy functional with respect to the \emph{Wasserstein} metric. 
The proof of the main result made no use of this structure, however, and this led us to ask,  \emph{Can we prove the same result using the structure of the Wasserstein gradient flow?}

This question is interesting for a number of reasons. The first is that the Wasserstein gradient flow is a \emph{natural} and \emph{physically meaningful} structure for this problem---we explain in Section~\ref{sec:ldp} what we mean by this. It can actually be argued that it is more natural than the linear structure that we used in the proof in~\cite{PeletierSavareVeneroni10}, and therefore it is also natural to ask whether this structure can be used.

The second reason is that the Wasserstein gradient-flow structure is known to arise in an impressively wide range of models and systems (e.g.~\cite{CarlenGangbo04,AmbrosioGigliSavare05,Savare07,BlanchetCalvezCarrillo08,MatthesMcCannSavare09,GianazzaSavareToscani09,Gigli10,CarrilloDiFrancescoFigalliLaurentSlepcev11}, just to name a few), and therefore any method that uses only the properties of this structure has the potential of application to a wide range of problems. Consequently, our approach here is to limit our use of information to those properties that follow directly from the gradient-flow structure. 

As a third reason, this work fits into a general endeavour to use
gradient-flow structures to pass to the limit in nonlinear
time-evolving systems (see
e.g.~\cite{SandierSerfaty04,Stefanelli08,MielkeRoubicekStefanelli08,MielkeStefanelli11,AmbrosioSavareZambotti09,Serfaty09TR}). 
The inherent convexity and lower-semicontinuity properties of this type of formulation provide handles for such limit passages that are similar to the well-known results for elliptic systems---as we show below.

\subsection{Kramers' problem}

The motion of a Brownian particle in a one-dimensional potential
landscape is described by the initial boundary-value problem (often called a Fokker-Planck or Smoluchovski equation~\cite[p.~8]{Risken84})
\begin{subequations}
\label{pb:main}
\begin{align}
\label{pb:main:eq}
&\partial_t \rho_\e = \tau_\e \partial_{\xi}\Bigl(\partial_\xi  \rho_\e + \frac1\e \rho_\e\partial_\xi H\Bigr),
&&t\geq 0,\ \xi\in \Xi:=[-1,1],\\
&\partial _\xi \rho_\e+\frac 1\e \rho_\e\partial_\xi H= 0, && t\geq 0, \ \xi=\pm 1.
\end{align}
\end{subequations}
The unknown function $\rho_\e$ is a time-dependent measure in
$\M(\Xi)$ (the space of finite, nonnegative, Borel measures
on the closed interval $\Xi=[-1,1]$), and this equation is to be interpreted in an appropriate weak form.

\begin{figure}[ht]
\centering
\noindent
\psfig{figure=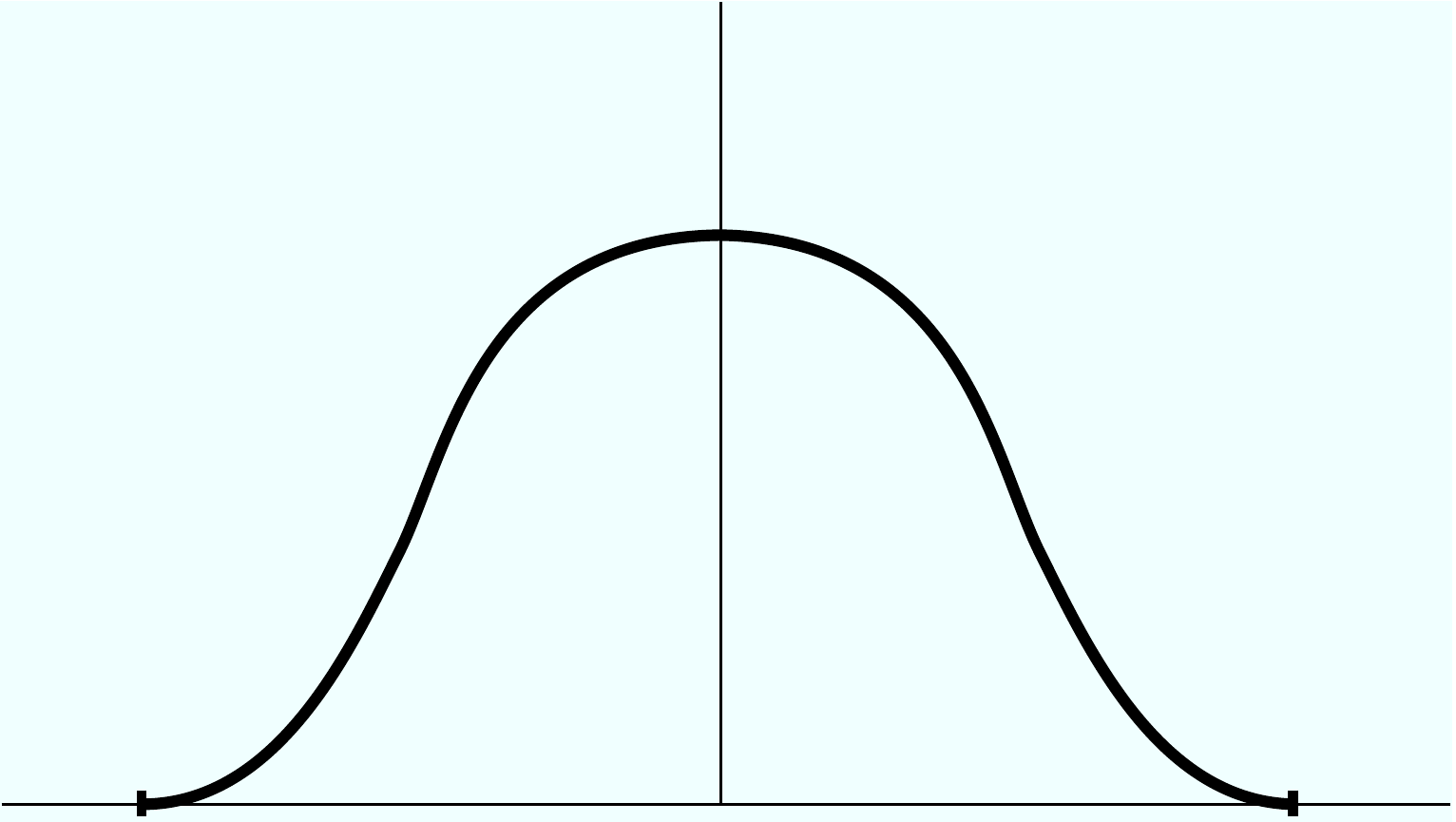,height=2cm}
\caption{A typical function $H$}
\label{fig:H}
\end{figure}

In this paper we take the potential energy $H$ to be a double-well potential, with wells in $\xi=\pm1$, and we follow the choice of~\cite{PeletierSavareVeneroni10} to truncate the domain at the wells, i.e. we take $\Xi=[-1,1]$ as the spatial domain (see Figure~\ref{fig:H}).
For definiteness we assume that $H$ is smooth, even, maximal at~$0$
with $H(0)=1$, and minimal at $\pm 1$ with $H(\pm1)=0$.
Each of these assumptions can be relaxed, but that is not the purpose of this paper.

In~\eqref{pb:main} two important constants appear. The potential $H$ is scaled by $1/\e$, which creates the situation of \emph{large activation energy}: the energy barrier separating the two wells is large in the limit $\e\to0$. As a consequence, the rate at which a particle passes from one well to the other is exponentially small as $\e\to0$; with the coefficient $\tau_\e$, Êwhich is defined in~\eqref{def:tau_e} below and which tends to infinity as $\e\to0$, we adapt the time scale to make the rate of transition asymptotically~$O(1)$.

\medskip

\begin{figure}[ht]
\centering
\noindent
\psfig{figure=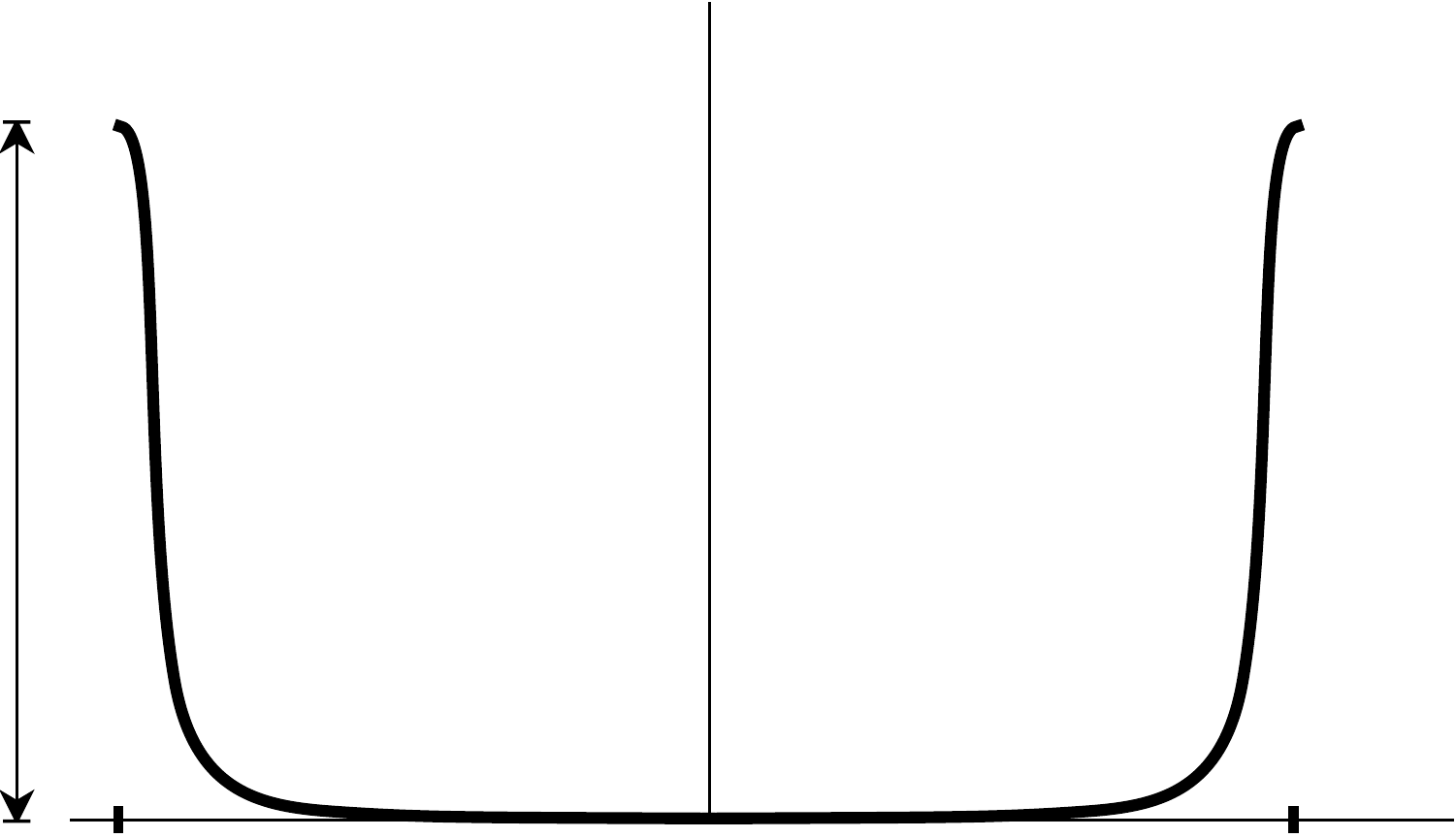,height=2.5cm}
\caption{The measure $\gamma_\e$, illustrated by plotting its Lebesgue density}
\label{fig:gamma-e}
\end{figure}

The asymptotically large `hump' of the potential $H/\e$ causes any solution of~\eqref{pb:main} to become singular in the limit $\e\to0$. This is well illustrated by the unique stationary solution  of unit mass,
\begin{equation}
\label{def:gamma_e}
\gamma_\e = Z_\e^{-1} \rme^{-H/\e}\L^1\restr{[-1,1]},
\end{equation}
where $Z_\e$ is a normalization constant and $\L^1$ is the one-dimensional Lebesgue measure (see Figure~\ref{fig:gamma-e}).
Since $H(\xi)>0$ for all $\xi\not=\pm1$, the measure 
 $\gamma_\e$ becomes strongly concentrated at the wells $\xi=\pm1 $ as $\e\to0$: 
\begin{equation}
  \label{conv:gamma_e}
  \gamma_\e\weakstarto \gamma_0=\frac 12\delta_{-1}+\frac 12\delta_1.
\end{equation}

In~\cite{PeletierSavareVeneroni10} we proved a number of results. The
first is that
the sequence $\rho_\e$ converges\footnote{The result of~\cite{PeletierSavareVeneroni10} uses a slightly different definition of $\tau_\e$, which is asympotically equivalent to the one this paper, \eqref{def:tau_e}.}, in the sense of measures, to a limit measure $\rho_0$, whose support is restricted to the two points $\xi=\pm1$:
\[
\rho_\e \weakstarto \rho_0 = \frac12 u_0^-\delta_{-1} + \frac12 u_0^+\delta_{1}.
\]
The densities $u_0^\pm:[0,T]\to\R$ of this limit measure $\rho_0$ satisfy the limit equation
\begin{subequations}
\label{pb:limit}
\begin{align}
\label{pb:limita}
\partial_t u_0^- &= k (u_0^+-u_0^-)\\
\partial_t u_0^+ &= k(u_0^- - u_0^+).
\end{align}
\end{subequations}
where the rate constant $k$ is given in terms of the potential function $H$ by
\begin{equation}
k=\frac1\pi \sqrt{|H''(0)|H''(1)}.\label{def:k}
\end{equation}
This limit system corresponds to the natural modelling of the monomolecular reaction $A\leftrightharpoons B$ at the continuum level.

\begin{figure}[h]
\begin{center}
\vskip4\jot
\noindent
{\psfig{height=2.5cm,figure=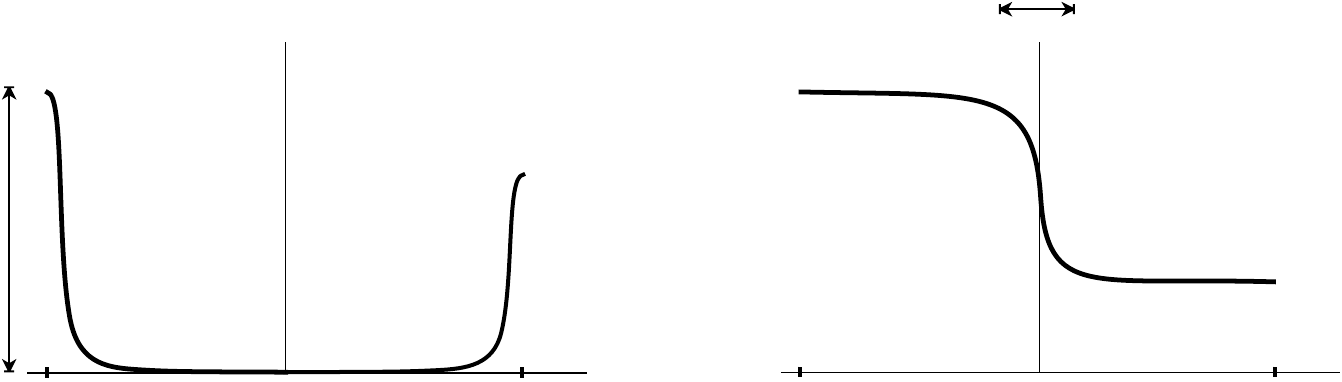}}
\end{center}
\caption{A comparison of $\rho_\e$ and the density $u_\e=\mathrm d\rho_\e/\mathrm d\gamma_\e$.}
\label{fig:ge-ue}
\end{figure}

A second result states a stronger form of convergence, and also highlights the role of the \emph{density}~$u_\e$ of the measure $\rho_\e$ with respect to $\gamma_\e$, i.e.
\begin{equation}
\label{def:ue}
u_\e = \frac{\mathrm d \rho_\e}{\mathrm d \gamma_\e}, 
\end{equation}
which satisfies the dual equation 
\begin{equation}
\label{pb:main_u}
\partial_t u_\e  = \tau_\e\Bigl(\partial_{\xi\xi} u_\e - \frac1\e \partial_\xi u_\e \partial_\xi  H\Bigr).
\end{equation}
Figure~\ref{fig:ge-ue} illustrates the relationship between $\rho_\e$ and $u_\e$. 
As it turns out, $u_\e$ is much better
behaved than $\rho_\e$ in the limit $\e\to0$: if the initial datum for
$u_\e$ is bounded above and below, then the same holds for $u_\e$
by the comparison principle, since constants are solutions of~\eqref{pb:main_u}.
In addition, $u_\e$ becomes locally constant away from $\xi=0$ (see part~\ref{th:compactness-b:locconstant} of Theorems~\ref{th:compactness-bounded} and~\ref{th:compactness} below).
This is reflected in a stronger form of convergence for $u_\e$, proved in~\cite{PeletierSavareVeneroni10}, which implies in particular that nonlinear functions of $u_\e$ also converge.

\bigskip

The aim of this paper is to derive similar convergence  statements by
different methods, specifically, by using only the
structure of the \emph{Wasserstein gradient flow}.
Before describing this structure for the specific case of~\eqref{pb:main}, we first recall the general structure of a gradient flow in a smooth and finite-dimensional setting.

%We describe this structure below, but first we introduce some functional-analytic notions that we need to characterize the limit behaviour of $\rho_\e$ in more detail. 

\subsection{Gradient flows {in a smooth Riemannian setting}}
  % and their
  %   variational characterization
\label{subsec:differentstructures}
Let
  us consider
  a smooth $d$-dimensional
  Riemannian manifold $\Z$, a $C^1$ energy functional 
  $\E:\Z\to\R$, and a quadratic dissipation potential $\psi$
  induced by the Riemannian metric on $\Z$.
  In local coordinates, we can identify $\Z$ (and the tangent space
  $\Tan z\Z$ at
  each point $z\in \Z$) with $\R^d$
  endowed with a smooth Riemannian tensor $G(z):\R^d(=\Tan z\Z)\to
  \R^d(=\CoTan z\Z)$ in the form $\psi(\dot
  z;z)=\frac 12 \langle G(z)\dot z,\dot z\rangle$.

The gradient flow of $\E$ in $\Z$ is
then given in the form 
\begin{equation}
  \label{eq:GradFlow}
  \dot z(t)=v(t)\in \Tan{z(t)}\Z,\quad\text{where}\quad
  v(t)=-\gradG \E(z(t))\quad\text{or}\quad
  G(z(t))v(t) = - \mathrm D \E(z(t))
  % \quad \text{or} \quad \dot z =
% -\nabla_{\!G} \E(z)
.
\end{equation}
Here and elsewhere in this paper 
we use overdots for time differentiation and $\mathrm D$ for the Fr\'ech\`et derivative of a function
(an element of $\CoTan z\Z$ in the Riemannian setting).
The  gradient $\nabla_{\!G}\E$ is defined as usual via the metric
as $z \mapsto G(z)^{-1} \mathrm D \E(z)$. It will sometimes  be easier to
use the dual dissipation potential $\psi^*$ given via the Legrendre
transform with respect to $\dot z$, namely
$\psi^*(\eta;z)=\frac12\langle \eta, G(z)^{-1} \eta \rangle$. Then
the gradient flow \eqref{eq:GradFlow} takes the form 
\begin{equation}
\label{eq:GradFlowPsi}
\dot z = \mathrm D\psi^*(-\mathrm D \E(z);z). 
\end{equation}
Here and below the derivatives $\mathrm D\psi$ and $\mathrm D\psi^*$
are only taken with respect to the first variable.

Solutions of \eqref{eq:GradFlowPsi} in a
  time interval $(a,b)$
  can be characterized as minimizers of the action functional
  \begin{equation}
    \label{eq:28}
    \cA(z;a,b):= \int_a^b \psi(\dot z+\gradG \E(z);z)\,\d t=
    \frac 12 \int_a^b \langle G(z) (\dot z+\gradG\E(z)),\dot
    z+\gradG\E(z)\rangle\,\d t,
  \end{equation}
  defined
  on $C^1$ curves with values in $\Z$.
  Expanding the integrand and observing that 
  \begin{displaymath}
    \langle G(z)\dot z,\gradG\E(z)\rangle=\langle \mathrm
    D\E(z),\dot z\rangle=\frac\d{\d t}\E(z),
  \end{displaymath}
  we see that $\cA$ has the structure
  \begin{equation}
    \label{eq:29}
    \begin{aligned}
      \cA(z;a,b)= \vphantom{A}&\E(z(b)) -\E(z(a)) + \cJ(z;a,b),\\
      \cJ(z;a,b):=\vphantom{A}&\int_a^b \Bigl[\psi\bigl(\dot z(t); z(t)\bigr) +
      \psi^*\bigl(-\mathrm D\E(z(t)); z(t)\bigr)\Bigr]\, \d t.
    \end{aligned}
  \end{equation}
  Note that for every curve $z$ we have
  $\cA(z;a,b)\ge0$, while $\cA(z;a,b)=0$ if and only if $z$ satisfies~\eqref{eq:GradFlow}.

  \subsection{Gradient flows in a metric setting}
  \label{subsec:metric-gradient-flows}
  The functionals $\cJ$ and $\cA$ can be generalized to infinite-dimensional 
  and non-smooth settings given by a
  space $\Z$ with  
  a lower semicontinuous (pseudo-, i.e.\ possibly taking the value
  $+\infty$) distance $\sfd:\Z\times \Z\to
  [0,+\infty]$.
  In such a space both tangent spaces and derivatives might not exist. Instead one can turn to two metric concepts,
  the metric slope $|\partial\E|$ of the functional $\E$ and the metric
  velocity $|\dot z|$ of a curve.
  The metric slope generalizes $(2 \psi^*(-\mathrm
  D\E(z);z) )^{1/2} $ and is defined by
  \begin{equation}
    \label{eq:38}
    |\partial\E|(z):=\limsup_{w\to
      z}\frac{(\E(z)-\E(w))_+}{\sfd(z,w)}. 
  \end{equation}
  Instead of defining a dissipation potential $\psi$ on the tangent
  space of an arbitrary point of $\Z$,
  one considers the class $AC(a,b;(\Z,\sfd))$ of absolutely continuous
  curves
  (with respect to the distance~$\sfd$)
  and their metric velocity
  \begin{equation}
    \label{eq:39}
    |\dot z|(t):=\lim_{h\to0}\frac {\sfd(z(t),z(t+h))}{|h|}\quad
    \text{if }{z\in AC(a,b;(\Z,\sfd))},
  \end{equation}
  which exists for a.e.\ $t\in (a,b)$~\cite[Th.~2.1.2]{AmbrosioGigliSavare05}.
  
  Using these concepts, the natural generalization of $\J$ in~\eqref{eq:29} is 
  \begin{equation}
    \label{eq:40}
    \J(z;a,b):=\int_0^T \Big[\frac 12|\dot z|^2(t)+\frac
    12|\partial\E|^2(z(t))\Big]\,\d t\quad\text{if }z\in AC(a,b;(\Z,\sfd)),
  \end{equation}
  trivially extended by $+\infty$ if $z$ is not absolutely continuous.
  Assuming that the slope is a \emph{strong upper gradient} for $\E$
  \cite[Ch.~2]{AmbrosioGigliSavare05},
  it is not difficult to prove that
  \begin{equation}
    \label{eq:41}
    \J(z;a,b)\ge \E(z(a))-\E(z(b))\quad\text{for every curve $z\in
      C([a,b];\Z)$ with $\E(z(a))<+\infty$.}
  \end{equation}
Comparing with the classical case outlined in Section~\ref{subsec:differentstructures} we deduce the following common structure: 
\begin{definition}
\label{def:GF}
Let $\Z$ be a topological space, $\E:\Z\to (-\infty,+\infty]$ be
a functional, and let $\cJ(\,\cdot\,;a,b)$ be a nonnegative (extended)
real functional
defined on $C([0,T];\Z)$ for all $0\leq a< b\leq T$, and satisfying 
\begin{equation}
  \label{cond:nonnegativeA}
  \E(z(b))+\J(z;a,b)\ge \E(z(a))\quad\text{for every
  }z\in C([a,b];\Z).
\end{equation}
Writing 
\[
\A(z) := \E(z(T)) - \E(z(0)) + \J(z;0,T),
\]
we define a curve $z\in C([0,T];\Z)$ to be a solution of the gradient flow system
$(\Z,\E,\cJ)$ if $\E(z(0))<\infty$ and $\A(z)=0$.
\end{definition}

This formulation of a gradient flow, in terms of the functional $\A$,
will be the basis for the rest of this paper. It clearly 
contains the classical case of a gradient system $(\Z,\E,\sfd)$, for
which $\J$ can be defined via \eqref{eq:38}--\eqref{eq:40}, and the metric-space outlined above, and it is sufficiently general to contain also the structure of the limiting problem (see Section~\ref{subsec:structureJ0}).
%However,
%while the system at positive $\e$ will be of this classical form, the limiting functional $\J_0$ in \eqref{def:J0} will have a different structure.  

\subsection{A first gradient-flow structure for~\eqref{pb:main}: the Hilbertian
  approach of \cite{PeletierSavareVeneroni10}}

We now turn to the specific case of this paper, equation~\eqref{pb:main}.
It is well known~\cite{Brezis73} that equation~\eqref{pb:main} in
the density formulation~\eqref{pb:main_u} is the 
gradient flow of the Dirichlet form 
\begin{equation}
  \label{eq:34}
  \E_\e^\text{lin}(u):=\frac {\tau_\e}2 \int_{-1}^1 
  \left|\partial_\xi u\right|^2 \,\dm{\gamma_\e}\xi,
\end{equation}
in the weighted Hilbert space $\Z^\text{lin}_\e=L^2(\Xi;\gamma_\e)$.
In this approach the quadratic dissipation potential (which is  also the
squared metric velocity) of a curve $u$ is 
\begin{equation}
\label{eq:GS_lin}
\begin{aligned}
\psi_\e^\text{lin}(\dot u;u)=\frac {1}2 \int_{-1}^1 
\dot u^2 \,\dm{\gamma_\e}\xi =\frac {1}2\|\dot u\|^2_{L^2(\Xi;\gamma_\e)},
\end{aligned}
\end{equation}
and does not depend on $u$, so that the resulting space has a flat
Hilbertian geometry.

The limit ODE \eqref{pb:limit} has a similar linear structure,
given by
\begin{equation}
\label{eq:GS_0_lin}
\begin{aligned}
&\Z_0^\text{lin}\;=\;L^2(\{-1,1\};\gamma_0)\;=\;\bigl{\{}\; (u_0^+,u_0^-)\in \R^2:
u_0^\pm=u_0(\pm 1)%  \;| \; u_0^+,u_0^-\geq 0,\ u_0^++u_0^-
% =\mass
\;\bigr{\}} ,\\
&\E_0^\text{lin}(u_0)\;=\;\frac k2 (u_0^+ {-} u_0^-)^2,\quad 
\psi_0^\text{lin}(\dot u_0; u_0)\;=\; \frac 14 |
\dot u^+_0|^2+ \frac 14|\dot u_0^-|^2\;=\;\frac 12\int_{\{-1,1\}} |\dot u_0|^2\,\dm{\gamma_0}\xi.
\end{aligned}
\end{equation}
The rigorous transition from \eqref{eq:GS_lin} to \eqref{eq:GS_0_lin} is
established in \cite{PeletierSavareVeneroni10} in a more general setting where 
diffusion in physical space is allowed as well. The 
analysis in \cite{PeletierSavareVeneroni10} 
depends in a crucial way on the
linearity of the problem.

\subsection{An alternative gradient-flow structure for~\eqref{pb:main}: the
  Wasserstein approach of \cite{JordanKinderlehrerOtto97}}
\label{subsec:intro_Wass}
As was discovered in the seminal work by
Otto and co-workers~\cite{JordanKinderlehrerOtto97,Otto01},  
equation \eqref{pb:main} has another relevant
gradient structure. 
It
relies on the interpretation of $\rho $ as a mass distribution
which is transported such as to reduce the free energy.

 In order to describe this point of view, we introduce
\begin{equation}
\label{eq:GS_Wass}\begin{aligned}
\Z^\text{meas}:= \M(\Xi),\qquad
 \E^\text{free}_\e(\rho):= \int_\Xi u \log u  \,
 \dm{\gamma_\e}\xi-\rho(\Xi)\log\rho(\Xi),\quad
 \text{where }u:=\frac{\d\rho}{\d\gamma_\e},
 % {\psi^*(\eta;u)=\int_{-1}^1 \frac {\tau_\e u}{2} |\partial_\xi
% \eta|^2\, \gamma_\e(d\xi). 
\end{aligned}
\end{equation}
with the convention that $\E^\text{free}(\rho)=+\infty$ if $\rho$ is
not absolutely continuous with respect to~$\gamma_\e$.
The space $\Z^\text{meas}$ is endowed with the usual weak-\textasteriskcentered \ convergence of measures (i.e.\
convergence in duality with continuous functions) and
can be metrized by the  \emph{$L^2$-Wasserstein distance}~$\sfd_W$.

This distance $\sfd_W$ admits two nice characterizations: the first one
involves optimal
transport
(see e.g.~\cite{Villani03,AmbrosioGigliSavare05}), while the second one is related to the dynamical interpretation discovered by
{Benamou and Brenier} \cite{BenamouBrenier00}
and is well adapted to the gradient-flow setting. 

In the latter point of view, we introduce the class
$\CE ab\Xi$ (\textbf{C}ontinuity \textbf Equation) given by couples $\rho\in C([a,b];\Z^{\text{meas}})$ and
$\nu\in \M((a,b)\times\Xi;\R)$ such that
\[
\partial_t \rho+\partial_\xi\nu=0\ \text{ in the sense of
    distributions in }\mathscr D'((a,b)\times\R).
\]
Here we trivially extend $\rho$ by zero outside of $\Xi$.
Often $\nu = \rho v$ for some Borel velocity field 
$v:(a,b)\times\Xi\to\R$, in which case the conditions above reduce to 
\begin{equation}
  \label{eq:35}
\begin{aligned}
& \int_a^b\!\! \int_\Xi |v(t,\xi)|\,\rho(t,\d\xi)\,\d t<+\infty \quad
 \text{ and }\\
& \partial_t \rho+\partial_\xi(\rho v)=0\ \text{ in the sense of
    distributions in }\mathscr D'((a,b)\times\R).
\end{aligned}
\end{equation}
%In the usual distributional formulation of the previous continuity
%equation we trivially extend $\rho$ and $v$ to $0$ outside $\Xi$.
For those couples $(\rho,\nu)\in \CE ab\Xi$ such that there exists such a velocity field $v$ with $\nu=\rho v$, the distance $\sfd_W$ can be defined in terms of $v$, by 
\begin{equation}
  \label{def:d_W}
  \sfd_W^2(\rho_0,\rho_1)= \min\biggl{\{}\int_0^1\!\! \int_{\bar D}
  \vert v(t, \xi)\vert^2\,\rho(t,\d \xi)\,\d t:(\rho,\rho v)\in\CE 01\Xi, \rho(0) = \rho_0,\rho(1) = \rho_1\biggr{\}},
\end{equation}
which illustrates how we can interpret $v$ as the `Wasserstein velocity' of the curve
$\rho$. 
Note how finiteness of $\sfd_W$ requires that $\nu \ll \rho$ and $\d\nu/\d\rho \in L^2(\rho)$, implying that
%Note the difference in powers of $v_t$ in~\eqref{def:d_W} and~\eqref{eq:35}, which shows how 
$\CE ab\Xi$ is a larger space than $AC([a,b];\sfd_W)$; indeed, our choice to work with $\CE ab\Xi$ stems from the fact that in the limit $\e=0$ the objects will still be elements of $\CE ab\Xi$, but no longer of $AC([a,b];\sfd_W)$. 
%The limiting measure~$\nu$ will no longer be absolutely continuous with respect to $\rho$, and therefore will not be of the form $\rho v$.

Recalling \eqref{eq:39}, it is natural to 
introduce the dissipation potential
\begin{equation}
  \label{eq:36}
  \psi^\text{Wass}_\e(v;\rho):=\frac {1}{2\tau_\e}\int_\Xi
  v(\xi)^2\,\dmm{\rho}\xi,\quad \text{for }\rho\in \M(\Xi),\ v\in L^2(\Xi;\rho).
\end{equation}
This expression suggests to interpret $L^2(\Xi;\rho(t,\cdot))$ as the
`Wasserstein tangent
space' at the measure~$\rho(t,\cdot)$, and in~\cite[Ch.~8]{AmbrosioGigliSavare05} this suggestion is made rigorous. 

The corresponding (squared) slope of $\E^\text{free}_\e$
\cite[\S 10.4.4]{AmbrosioGigliSavare05} defined by
\eqref{eq:38} is the \emph{Fisher information}
\begin{equation*}
  \label{eq:35bis}
  |\partial\E^\text{free}(\rho)|^2:=\int_\Xi \Big|\frac{\partial_\xi
    u}u\Big|^2\,\dm\rho\xi=
  4\int_\Xi \big|\partial_\xi \sqrt u\big|^2\,\dm{\gamma_\e}\xi,\quad\text{if
  }u=\frac{\d\rho}{\d\gamma_\e}\text{ with } \sqrt u\in W^{1,2}(-1,1).
\end{equation*}
This corresponds to the choice of the dual dissipation potential
\begin{equation*}
  \label{eq:31}
  (\psi^\text{Wass}_\e)^*(\eta;\rho):=\frac {\tau_\e}2\int_\Xi
  |\eta(\xi)|^2\,\dmm\rho\xi
\end{equation*}
and of the `Wasserstein gradient' $\nabla_{\!W}\E^\text{free}_\e$ of
the entropy given (at least formally) by
\begin{equation*}
  \label{eq:42}
  \nabla_{\!W}\E^\text{free}_\e(\rho):=\frac{\partial_\xi
    u}u=\partial_\xi\log u,\quad u=\frac {\d\rho}{\d\gamma_\e}.
\end{equation*}
This construction is equivalent to~\eqref{pb:main}: in fact, at least for
smooth densities, 
\begin{equation*}
  \label{eq:44}
  \partial_\xi\rho+\frac 1\e\rho\partial_\xi H=\rho\,\partial_\xi \log\Big(\frac {\d\rho}{\d\gamma_\e}\Big)
\end{equation*}
so that \eqref{pb:main} has the gradient flow structure
\eqref{eq:GradFlow} in the Wasserstein sense:
\begin{equation*}
  \label{eq:43}
  \partial_t \rho+\partial_\xi(\rho\, v)=0,\qquad
  v=- \tau_\e\nabla_{\!W}\E^\text{free}_\e(\rho) ,
\end{equation*}

Motivated by these remarks, we introduce the 
functional $\J_\e^\text{Wass}$, 
\begin{multline}
  \label{def:JWass}
  \J_\e^\text{Wass}(\rho;a,b):=\int_a^b\Big(\int_\Xi
  \frac1{2\tau_\e}v^2\,\rho(t,\d\xi)+
  2\tau_\e%\frac {\tau_\e}2
  \int_\Xi \big|\partial_\xi
  \sqrt{u}\big|^2\,\d\gamma_\e\Big)\,\d
  t,\\
  \quad\text{if }(\rho,\rho v)\in \CE ab\Xi \text{ and }\rho=u\gamma_\e,
\end{multline}
and the corresponding Wasserstein action functional
\begin{equation}
  \label{def:AWass}
  \A_\e^\text{Wass}(\rho;a,b):=\E_\e^\text{free}(\rho(a))-\E_\e^\text{free}(\rho(b))+\J_\e^\text{Wass}(\rho;a,b),
\end{equation}
which satisfies the admissibility condition
\eqref{cond:nonnegativeA}.
In analogy to Definition \ref{def:GF},
a Wasserstein solution $\rho$ of \eqref{pb:main} in the time interval
$[0,T]$ is a curve in
$\M(\Xi)$ with $\E_\e^\text{free}(\rho(0))<+\infty$ and
$\A_\e(\rho;0,T)=0$.

\subsection{Our main results}
\label{subsec:results}

In this work we prove various results on the connection between the Wasserstein gradient
structure $(\Z^\text{meas},\E_\e^\text{free},\J_\e^\text{Wass})$ 
and a gradient structure 
$(\Z_0^\text{meas},\E_0^\text{free},\J_0)$ for the limit
system~\eqref{pb:limit}.
As described above, the motivating question is whether we may pass to the limit in the gradient-flow equation $\A^\text{Wass}_\e(\rho_\e)=0$.
This question falls apart
into two sub-questions: 
\begin{enumerate}
\item \emph{Compactness:} Do solutions of
  $\cA^\text{Wass}_\e(\rho_\e;0,T)=0$ with uniformly bounded initial
  entropy $\E_\e^\text{free}(\rho_\e(0))$
  have beneficial
 compactness properties, allowing us to extract a subsequence that
 converges in a suitable topology, say $\sigma$?
\item \emph{Liminf inequality:} Is there a limit functional $\J_0$
  such that
  \begin{equation*}
%    \label{eq:50}
    \rho_\e\stackrel\sigma\longrightarrow\rho_0\quad\Longrightarrow\quad 
    \liminf_{\e\to0}\J_\e^\text{Wass}(\rho_\e;a,b)\ge \J_0(\rho_0;a,b)?
  \end{equation*}
  And if so, does it satisfy the admissibility condition \eqref{cond:nonnegativeA}, i.e.
  \begin{equation*}
%    \label{def:A0}
%    \cA_0(\rho;a,b):=
\E_0^\text{free}(\rho(a))-\E_0^\text{free}(\rho(b))+\J_0(\rho;a,b)\ge
    0
  \end{equation*}
  for every $0\leq a < b\leq T$ and $\rho\in C([a,b];\Z_0^\text{meas})$ with
  $\E_0^\text{free}(\rho(a))<+\infty$?
  \end{enumerate}

Our answer to these questions is indeed affirmative. 
Question 1 is answered by Theorem~\ref{th:compactness}, which establishes that any sequence
$\rho_\e$ such that $\E_\e^\text{free}(\rho_{\e}(0))$ and $\J_\e^\text{Wass}(\rho_\e;0,T)$
are bounded, is compact in several topologies. These boundedness assumptions are natural, and only use information associated with the gradient-flow structure. 

%The bound on $\E_\e(u_\e(0))$ is a natural requirement, and corresponds to taking a reasonable sequence of initial data. The $L^\infty$ bound is less natural, and we would prefer to remove it; technical issues prevent that at the moment. 
Question 2 is addressed by Theorems~\ref{th:inner_lower_bound} and~\ref{th:recovery_sequence}, which 
characterize the limit of the functionals $\J_\e^\text{Wass}(\cdot;a,b)$ in terms of
Gamma-convergence.
{If we denote by $\sigma$ the topology mentioned above},
then this convergence is characterized by
the existence of functionals $\J_0$ and $\E_0^\text{free}$ satisfying the two properties
\begin{enumerate}
\item Lower bound: for each family of curves
  $\rho_\e\stackrel\sigma\longrightarrow \rho_0$
  with $\sup_\e \E^\text{free}_\e(\rho_\e(a))<+\infty$, we have
  
\[
\J_0(\rho_0;a,b)\leq \liminf_{\e\to0}
\J_\e^\text{Wass}(\rho_\e;a,b)\quad\text{and}\quad
\E^\text{free}_0(\rho_0(b))\le \liminf_{\e\to0}\E^\text{free}_\e(\rho_\e(b)).
\]

\item Recovery sequence: for each
   $\rho_0\in C([a,b];\Z^\text{meas}_0)$
  with $\E^\text{free}_0(\rho_0(a))<+\infty, \
  \J_0(\rho_0;a,b)<+\infty$
  there exists a sequence
  $\rho_\e\in C([a,b];\Z^\text{meas})$
  such that $\rho_\e\stackrel\sigma\longrightarrow \rho_0$ and
\[
\J_0(\rho_0;a,b) = \lim_{\e\to0} \J_\e^\text{Wass}(\rho_\e;a,b),\quad
\E^\text{free}_0(\rho_0(b))= \lim_{\e\to0}\E^\text{free}_\e(\rho_\e(b)).
\]
\end{enumerate}

The limit structure $(\Z_0^\text{meas},\E_0^\text{free},\J_0)$ consists of measures $\rho$ that are absolutely continuous with respect
to $\gamma_0$ and thus supported in $\{-1,1\}$:
\begin{equation*}
\rho=\frac 12 u^-\delta_{-1}+\frac 12 u^+\delta_1 
\qquad\text{for some }u^\pm\geq0.
%\label{eq:48}
\end{equation*}
The space $\Z_0^\text{meas}$ and the energy $\E_0^\text{free}$ are natural limits of the corresponding objects as $\e\to0$:
\begin{equation}
\label{eq:GS_0_Wass}\begin{aligned}
&\Z_0^\text{meas}= \Big\{\rho\in \M(\Xi):\supp(\rho)\subset
\{-1,1\}\Big\}
% =\biggl\{\; \rho =\frac12 u_0^-\delta_{-1} + 
% \frac12 u_0^+ \delta_1 \in \M([-1,1])\;\biggr\}
\subset \Z^\text{meas},\\
& \E_0^\text{free}(\rho)= \int_{\{-1,1\}} \frac
{\d\rho}{\d\gamma_0}\log\Big(\frac
{\d\rho}{\d\gamma_0}\Big)\,\d\gamma_0-
\mass\log \mass=
\frac12\Big( u^+ \log u^+ + u^- \log
u^- \Big)-\mass\log\mass\\
&\text{where}\quad
u^\pm=\frac{\d\rho}{\d\gamma_0}(\pm 1),\quad
\mass=\rho(\Xi)=\frac 12(u^++u^-).
%  \quad  {\psi_0^\text{Wass}}^*(\eta_0^\pm; u_0^\pm)=
% \frac k2 \,\frac{u_0^+-u_0^-}{\log u_0^+- \log u_0^-}\, 
% | \eta_0^+- \eta_0^- |^2 ,  
\end{aligned}
\end{equation}
% where the second fraction in the definition of ${\psi_0^\text{Wass}}^*
% $ has the continuous extension to $u_0^+$ for the case
% $u_0^+=u_0^-.$
This limit energy $\E_0^\text{free}$ is the Gamma-limit of $\E_\e^\text{free}$~\cite{AmbrosioSavareZambotti09}. 
%This gradient structure has a generalization to much
%more general reaction-diffusion systems, see~\cite{Mielke10TR}. 
However, the limit functional $\J_0(\cdot;a,b)$, does not have the same duality structure as~\eqref{eq:29}, and we discuss this next.

\subsection{The structure of $\J_0$}
\label{subsec:structureJ0}
In fact, since the limit problem is characterized by measures $\rho(t)$ concentrated
at $\xi=\pm 1$,  no effective mass transport is possible between $\xi=-1$ and $\xi=1$.
Assume for instance the case when $\rho$ is sufficiently smooth, i.e. $\rho(t)=\frac 12\sum u^\pm(t)\delta_{\pm
  1}=u\gamma_0$
for a couple $u^\pm\in C^1((a,b);\R)$. Then the distributional
time derivative
of $\rho$ is $\partial_t \rho=\frac 12\sum \dot u^\pm \delta_{\pm
  1}$ and any signed measure $\nu$ supported in $\Xi\times [a,b]$ and solving the
continuity equation
\begin{equation}
  \label{eq:68}
  \partial_t \rho+\partial_\xi\nu=0\quad\text{in }\mathscr D'(\R\times (a,b))
\end{equation}
cannot be absolutely continuous with respect to $\rho$
and therefore cannot admit the decomposition $\nu=\rho v$ for some
$v\in L^2(-1,1;\rho)$
(except for the trivial case $\dot u^\pm\equiv 0$).

Recalling that the total mass $\mass=\frac12(u^-+u^+)$ is
conserved
and therefore $\dot u^-=-\dot u^+$, equation~\eqref{eq:68} has the unique solution
\begin{equation}
  \label{eq:61}
  \nu=w\L^2\restr{(-1,1)\times (a,b)},\quad
  w(\xi,t)=\frac 12 \dot u^+(t)\quad \text{for }\xi\in (-1,1),\ t\in (a,b),
\end{equation}
trivially extended to $0$ outside $[-1,1]$.

As we show below, $\J_0(\rho;a,b)$ has the form
\begin{align}
  \label{def:J0}
    \J_0&(\rho;a,b):=\int_0^T \Md{w(t)}{u^\pm(t)}\,\d t\\&\text{if
    }\rho(t)=\frac 12 \sum u^\pm(t)\delta_{\pm 1},\quad u^\pm \in
    AC(a,b;\R),\quad \text{with } w=\frac 12\dot u^- = -\frac12 \dot u^+,
    \notag
\end{align}
where the function $M:\R\times [0,\infty)^2\to [0,+\infty]$ is given by
\begin{equation}
\label{def:M}
\Md w{u^\pm} := \inf \left\{ \,\int_{-\kappa}^\kappa \frac{w^2}{2u(s)}  {+}
  \frac{u'(s)^2}{2u(s)}\: \d s\ :\ 
  u\in H^1(-\kappa,\kappa), \ u(\pm \kappa) = u^\pm \,\right\},\quad
\kappa:=\frac 1 k.
\end{equation}

This functional $\J_0$ satisfies the admissibility criterion~\eqref{cond:nonnegativeA}. 
Indeed,  along any admissible curve $\rho(t)$,
\begin{equation*}
  \label{eq:69}
  \begin{aligned}
    \frac \d{\d t}\E_0^\text{free}(\rho(t))&= \frac \d{\d t}\frac 12
    \sum_\pm u^\pm(t)\log u^\pm(t)= \frac 12\sum_\pm (\log
    u^\pm(t)+1)\dot u^\pm(t)\\
    &= \frac 12[\log u^+(t)-\log u^-(t)]\,\dot
    u^+(t),
  \end{aligned}
\end{equation*}
so that the admissibility condition 
\eqref{cond:nonnegativeA} is equivalent to
\begin{equation}
  \label{eq:70}
  (\log u^+-\log u^-)w\le \Md w{u^\pm}\quad\text{for every }w\in
  \R,\ u^\pm>0.
\end{equation}
In Theorem~\ref{th:M.est} we prove this inequality, implying that the limiting action $\A_0$, 
\begin{equation}
\label{def:A0}
\A_0(\rho) := \E_0^{\text{free}}(\rho(T)) - \E_0^{\text{free}}(\rho(0)) 
  + \J_0(\rho),
\end{equation}
satisfies $\A_0(\rho)\geq0$ for all $\rho$.

It is now natural to ask which curves $\rho$ satisfy the equation $\A_0(\rho;a,b)=0$. This equation implies equality in~\eqref{eq:70}, which suggests defining the `contact set'~\cite{MielkeRossiSavare09}
\[
\mathcal C:= \{ \; (u^\pm,w) \;|\; \Md w{u^\pm} + (\log u^+-\log
u^-)w=0
\;\}.
\]
A consequence of a second inequality proved in Theorem~\ref{th:M.est} is 
\[
(u^\pm,w)\in \mathcal C \quad \Longleftrightarrow \quad 
w = \frac k2(u^+-u^-).
\]
This implies that any $\rho$ satisfying $\A_0(\rho;a,b)=0$ also solves the limiting equation \eqref{pb:limit}. 

\subsection{Recovering a gradient flow}
Finally, one might ask whether it is possible to find a `true'
gradient structure, i.e. an alternative functional $\A_0$ that does
have the dual structure as in~\eqref{eq:29} or~\eqref{eq:40}. For this we need to find a dissipation potential
$\psi_0(w; u^\pm)$ such that the associated contact set is equal to
$\mathcal C$, i.e. such that %are the "??" below meant to be there?
\[
\mathcal C  =  \{\; (w,u^\pm)\;:\; \psi_0(w;u^\pm) +\psi^*_0(-\mathrm
D \E_0(u^\pm);u^\pm) + \langle w,\mathrm D \E_0(u^\pm) \rangle =0\;\}. 
\]
Using the two-sided estimate of Theorem~\ref{th:M.est} for $M$ we
find that a natural choice for $\psi_0$ is 
\[
 \psi_0(w;u^\pm)= \frac2k\frac{\log u^+ - \log u^-}{u^+-u^-}w^2
\]
which gives the desired result \eqref{pb:limit}.  

\subsection{The variational approach: basic tools and main ideas.}
As we mentioned before, in the present paper we adopted a
metric-variational approach to extract crucial information from the
particular structure of the equation \eqref{pb:main}. 
This point of view has become quite popular and arises from
\begin{itemize}
\item the combination of general metric concepts (briefly recalled
  in
  Section \ref{subsec:metric-gradient-flows}),
\item basic measure-theoretic tools, optimal transportation, and
  entropy-dissipation techniques (\`a la Otto \cite{Otto01},
  see Section \ref{subsec:intro_Wass}),
\item a careful use of the continuity equation \eqref{eq:35} and
  of the Benamou-Brenier dynamical point of
  view~\cite{BenamouBrenier00}, 
\item standard $\Gamma$-convergence methods.
\end{itemize}
On the other hand, the problem exhibits many non-standard features,
which we have addressed with new ideas and techniques, which could
hopefully be
useful in other situations.

\begin{itemize}
\item One of the main points here is that the limit procedure mixes in an unusual way
  the two contributions to $\J_\e$, namely that due to time change ($\psi(\dot z)$ in the notation of Section~\ref{subsec:differentstructures}) and that due to the entropy slope ($\psi^*(-\mathrm D \E)$). Therefore one cannot
  canonically separate these two contributions in the structure of $\J_0$. This fact has
  another interesting consequence: in the present setting it is not
  possible to investigate separately the limit behaviour of the
  distance and of the functional using $\Gamma$-convergence tools (as
  in the well-behaved gradient flows considered by
  \cite{Savare07,AmbrosioSavareZambotti09,Serfaty09TR,PeletierSavareVeneroni10}).
  Conversely, the geometry perturbed by the sublevels and by the
  slopes of the varying entropy functionals $\E_\e^{\rm free}$ induces
  a new kind of evolution in the limit, which can solely be captured
  by considering the asymptotics of the whole space-time functionals
  $\J_\e$.\\
  This singular behaviour motivated our general Definition
  \ref{def:GF} of a gradient flow system $(\A,\E,\J)$ and the idea to
  focus on basic structural dissipation inequalities along solutions
  of the continuity equation which could be preserved in the limit.

\item A second crucial point is the rescaling strategy (Section
  \ref{sec:rescaling}) which allows us to resolve the singular
  behaviour of the functionals along sharp transitions.  Since the
  continuity equation has good invariance properties, we can combine
  the information coming both from the original and from the rescaled
  formulation to construct the limit.
\item As a byproduct,
  we recover a (quite surprising, at first sight) common structure
  uniting \eqref{pb:main} and \eqref{pb:limit} thanks to this
  measure-theoretic interpretation: this is not completely trivial, if
  one takes into account that the limit problem is a system of
  ordinary differential equations.

\item  This fact is strongly related to the new variational formulation of
reaction-diffusion systems, see \cite{Mielke11,GlitzkyMielke11}.
In fact, it turns out that many reaction systems (even with nonlinear reactions
of mass-action type) can be written as a gradient system with a well-chosen
dissipation functional $\psi^*_0$. While for linear reaction systems such as here
(or more generally Markov systems) there are several variational gradient
structures, this new formulation appears to be the only one that is compatible with
diffusion processes.

\item A final comment should be devoted to the connections with 
  stochastic particle systems, which we will explain in more detail in
  Section \ref{sec:ldp}.
\end{itemize}
\subsection{Structure of the paper}

In Section~\ref{sec:rescaling} we introduce a rescaling of space that
desingularizes one of the terms in $\J_\e$. This rescaling allows us
to prove, in Section~\ref{sec:compactness}, the compactness of a
sequence $\rho_\e$ with bounded initial energy $\E_\e(\rho(0))$ and
bounded $\J_\e(\rho_\e)$ in a number of
topologies. Sections~\ref{sec:lowerbound} and~\ref{sec:upperbound}
give the two parts of the Gamma-convergence result, the lower bound
and the recovery sequence. Before constructing the recovery
sequence we investigate in Section~\ref{sec:M} the function $M$ in
some detail. These are the central mathematical results of the paper.

In Section~\ref{sec:ldp} we place the results of this paper in the
context of large-deviation principles for systems of Brownian
particles, and comment on the various connections. In   Section~\ref{sec:discussion} we discuss various aspects of the results and their proof and comment on possible generalizations. Finally, in Section~\ref{sec:HerrmannNiethammer} we draw parallels between this work and an independent study of the same question by Herrmann and Niethammer~\cite{HerrmannNiethammer11}.

\bigskip
{\parindent0pt
\begin{minipage}{\textwidth}
\centering{\bf Summary of notation}
\ \\[3mm]
\begin{small}
\begin{tabular}{lll}
$\longweakstarto$        & weak convergence in duality with continuous functions &\\
$\CE \cdot\cdot\cdot$ & pairs $(\rho,\nu)$ satisfying the continuity equation&\eqref{eq:35}\\
$\sfd_W^2(\cdot,\cdot)$ & Wasserstein distance of order 2  &\eqref{def:d_W}\\
$\E_\e$, $\J_\e$, and $\A_\e$ & (Section~\ref{sec:introduction}) general energy, dissipation fuctional, and action \\
$\E_\e$, $\J_\e$, and $\A_\e$ & (Sections~\ref{sec:rescaling}--\ref{sec:HerrmannNiethammer}) Wasserstein energy, dissipation fuctional, and action,\\
& i.e the same as $\E_\e^{\mathrm{free}}$, $\J_\e^{\mathrm{Wass}}$, and $\A_\e^{\mathrm{Wass}}$  & \eqref{eq:GS_Wass}, \eqref{def:JWass}, \eqref{def:AWass}\\
$\hE_\e$, $\hat\J_\e$, and $\hat \A_\e$ & $\E_\e$, $\J_\e$, and $\A_\e$ written in terms of $\hat \rho$ & \eqref{def:hatE}, \eqref{def:hatJ}, \eqref{eq:47bis} \\
$\E_0$, $\J_0$, and $\A_0$  & limit energy, dissipation, and action &\eqref{eq:GS_0_Wass}, \eqref{def:J0}, \eqref{def:A0}\\
$\gamma_\e$, $\hge$  & invariant measure ($\gamma_\e$) and its push-forward under $\hat s_\e$  &\eqref{def:gamma_e} and \eqref{def:hge}\\
$\gdens$, $\hgdens$ & Lebesgue densities of $\gamma_\e$ and $\hge$ &\eqref{eq:51}, \eqref{def:hgdens}\\
$H$                  & `enthalpy function', potential for the Brownian particle & page \pageref{fig:H}\\
$k=1/\kappa$         & reaction parameter & \eqref{def:k}\\
$\Md \cdot\cdot$                  & argument of the integral in $\J_0$ & \eqref{def:M}\\
$\hat s_\e$			& transform from $\xi$ to $s$, inverse of $\hat \xi_\e$ & \eqref{diffeo}\\
$\tau_\e$			& time rescaling & \eqref{def:tau_e}\\
$u_\e$               & density $\d \rho_\e/\d \gamma_\e$  & \eqref{def:ue}\\
$\hat u_\e$          & transform of $u_\e$, $\hat u_\e = u_\e \circ \xi_\e$ &\eqref{def:hatu}\\
$w_\e $				& Lebesgue density of $v_\e \rho_\e$ & \eqref{def:w_e} \\
$\hat w_\e$ 			& Lebesgue density of $\hat v_\e \hat \rho_\e$ &\eqref{eq:43bis}\\
$\hat \xi_\e$		& transform from $s$ to $\xi$, inverse of $\hat s_\e$ & \eqref{diffeo}\\
$Z_\e$				& normalization constant of $\gamma_\e$ & \eqref{def:gamma_e}
\end{tabular}
\end{small}
\end{minipage}
}

\section{Rescaling}
\label{sec:rescaling}

\subsection{Definitions}

From here on we write $\E_\e$, $\J_\e$, and $\A_\e$ for  $\E_\e^{\mathrm{free}}$, $\J_\e^{\mathrm{Wass}}$, and $\A_\e^{\mathrm{Wass}}$, since we will only be using the Wasserstein framework. Since for most of the discussion the interval $(a,b)$ will be fixed to $(0,T)$, we will also write $\J(\rho)$ for $\J(\rho;0,T)$ etcetera.

A central step in the analysis of this paper is a rescaling of the
domain which stretches the region around $\xi=0$. This converts the
functions $u_\e$, which have steep gradients around $\xi=0$ (see
Figure~\ref{fig:ge-ue}), into functions $\hue$ of the new variable $s$
that will have a more regular behaviour.

We call $\gdens$ the Lebesgue density of $\gamma_\e$, namely 
\begin{equation}
  \label{eq:51}
  \gdens(\xi):=Z_\e^{-1}\rme^{-H(\xi)/\e}, \text{ and } 
  \text{we set }   \kappa:=\frac 1 k=
  \frac{\sqrt{|H''(0)|\,H''(1)}}{\pi} \quad \text{(as in
    \eqref{def:k})} .
\end{equation}
We now make the choice of $\tau_\e$ precise:
\begin{equation}
\label{def:tau_e}
\tau_\e := \frac 1{2\kappa}\int_{-1}^1 \frac{\d\xi}{\gdens(\xi)}.
\end{equation}
An application of Watson's Lemma gives the asymptotic estimate
\[
\left.
   \tau_\e%\right.
%\left.
\middle/%
\e \rme^{1/\e} \right. \;\;\stackrel{\e\to0}\longrightarrow \;\;1.%\frac{2\pi}{\GG{\sqrt{|H''(0)|\,H''(1)}}}.
\]

Using that $\gdens$ is even (since $H$ is even), 
we introduce the smooth increasing diffeomorphism
$\xi\mapsto s= \hat s_\e(\xi)$,
\begin{equation}
\label{diffeo}
\hat s_\epsilon: [-1,1]\rightarrow [-\kappa,\kappa],\quad \hat s_\epsilon(\xi):= 
\int_0^\xi \frac{1}{\tau_\epsilon \gdens(\eta)}\,\d\eta, \
\quad \text{with inverse} \quad \hat \xi_\epsilon:= \hat
s_\epsilon^{-1} : [-\kappa,\kappa]
\rightarrow [-1,1] .
\end{equation}
Note that $\hat\xi_\e$ satisfies
\begin{equation}
  \label{eq:52}
  \frac \d{\d s}\hat\xi_\e(s)=\tau_\e \gdens(\hat\xi_\e(s)),\quad
  \hat \xi_\e(-\kappa)=-1 , \ \text{ and } \hat
  \xi_\e(\kappa)=1.
\end{equation}

With this change of variables we call $\rmS:=[-\kappa,\kappa]$
the domain of the variable $s$
and we set
\begin{equation}
  \label{def:hge}
  \hat\gamma_\e:=(\hat s_\e)_\#\gamma_\e,\quad
  \hat \rho_\e:=(\hat s_\e)_\#\rho_\e.
\end{equation}
Observe that the Lebesgue density $\hgdens$ of $\hge$
satisfies
\begin{equation}
\label{def:hgdens}
\hgdens(\hat s_\e(\xi)) \frac{\d}{\d\xi} \hat s_\e(\xi) = \gdens(\xi)
\quad\text{so   that}\quad
\hgdens(s)=\tau_\e \gdens^2(\hat \xi_\e(s)),
\end{equation}
and the transformed measure $\hat\rho_\e$ satisfies
\begin{equation}
\label{def:hatu}
\hat \rho_\e=\hat u_\e\hat \gamma_\e,\quad
\hat u_\epsilon := u_\epsilon\circ \hat \xi_\e.
\end{equation}
In particular, we can easily transport the entropy functional to the
new setting
\begin{equation}
\label{def:hatE}
\hE_\e (\hat \rho):= \int_\rmS \hat u(s)\log\hat u(s)\,
\dmm{\hat\gamma_\e }s,\quad
\text{so that }
\quad
\hE_\e(\hat\rho)=\E_\e(\rho)\quad\text{if }
\hat \rho=(\hat s_\e)_\#\rho.
\end{equation}
If $(\rho,\rho v)\in \CE ab\Xi$ then the couple $(\hat \rho,\hat \rho\hat v)$ with
$\hat v(s)=v(\hat \xi(s))/(\tau_\e\gdens(\hat \xi_\e(s)))$
belongs to $\CE ab\rmS$ and satisfies the continuity equation
\begin{equation}
  \label{eq:42bis}
  \partial_t \hat \rho+\partial_s (\hat \rho\,\hat v)=0\quad\text{in
  }\mathscr D'((0,T)\times\R);
\end{equation}
in fact, since 
$v(\hat \xi(s))=\hat v(s)\hat \xi_\e'(s)$, for every
$\hat\phi=\phi\circ\hat \xi\in C^\infty_{\rm c}((0,T)\times[0,1])$ we have
\begin{align*}
  0&=  \int_0^T\!\!\int_{-1}^1 \Big(\partial_t \phi+v\,\partial_\xi \phi 
  \Big)\,\dmm{\rho_t}\xi\,\d t
  =\int_0^T\!\!\int_\rmS \Big((\partial_t \phi)\circ\hat \xi_\e+
  (v\circ\hat\xi_\e)\,(\partial_\xi \phi) \circ\hat\xi_\e\Big)\,\dmm{\hat\rho_t}s\,\d t\\&=
  \int_0^T\!\!\int_\rmS \Big(\partial_t \hat\phi
    +\hat \xi_\e'\,\hat v\, (\partial_\xi
  \phi)\circ\hat \xi_\e   \Big)\,\dmm{\hat\rho_t}s\,\d t
=\int_0^T\!\!\int_\rmS \Big(\partial_t \hat\phi
  + \hat v \, \partial_s
  \hat \phi
   \Big)\,\dmm{\hat\rho_t}s\,\d t.
\end{align*}
Setting $\hat w:=\hat v\,\hat u \,\hgdens$ we also have
\begin{equation}
  \label{eq:43bis}
  \frac1{2\tau_\e}\int_\Xi v^2\,\d\rho=
  \frac1{2\tau_\e}\int_{\rmS} \hat v^2\,\tau^2_\e \gdens^2(\hat
  \xi_\e)\,\d\hat \rho=
  \frac1{2}\int_{\rmS} \hat v^2\,\hgdens^2\hat
  u\,\d s=
  \frac1{2}\int_{\rmS} \frac{\hat w^2}{\hat
  u}\,\d s .
\end{equation}
Since
\begin{equation*}
  \label{eq:44bis}
  \partial_s \sqrt{\hat u}=
  (\partial_\xi\sqrt u\circ \hat \xi_\e)\,\hat \xi_\e'=
  \tau_\e \,(\partial_\xi\sqrt u\circ \hat \xi_\e)\;\gdens\circ \hat \xi_\e
\end{equation*}
we also get
\begin{equation}
  \label{eq:45bis}
  2\tau_\e\int_\Xi \Big|\partial_\xi\sqrt u\Big|^2\,\d\gamma_\e=
  2\int_{\rmS} \Big|\partial_s \sqrt{\hat
    u}\Big|^2\frac{1}{\tau_\e \gdens^2(\hat\xi_\e)}\,\d\hat\gamma_\e=
  2\int_{\rmS} \Big|\partial_s \sqrt{\hat
    u}\Big|^2\,\d s.
\end{equation}
Combining \eqref{eq:42bis}, \eqref{eq:43bis}, and \eqref{eq:45bis}, we now
define the functional
\begin{equation}
  \label{def:hatJ}
  \hat\J_\e(\hat\rho;0,T):=\int_0^T\Big(\frac 12 \int_{\rmS}
  \frac{\hat w^2}{\hat
  u}\,\d s+2 \int_{\rmS} \Big|\partial_s \sqrt{\hat
    u}\Big|^2\,\d s\Big)\,\d t,\quad
  \hat u=\frac{\d\hat \rho}{\d\hat\gamma_\e},\quad
  \hat w=\hat u\,\hat v\,\hgdens,
\end{equation}
and $(\hat \rho,\hat\rho \hat v)=
(\hat \rho,\hat w\L^2)\in \CE ab\rmS$.
This calculation shows that
\begin{equation}
\begin{split}
  \label{eq:47bis}
  &\J_\e(\rho;0,T)=
  \hat \J_\e(\hat \rho;0,T),\quad\text{and}\\
  &\A_\e(\rho;0,T)=
  \hat A_\e(\hat \rho;0,T)=
  \hE_\e(\hat \rho(b)) - \hE_\e(\hat \rho(a)) +\hat \J_\e(\hat \rho;0,T).
\end{split}
\end{equation}

\begin{figure}[h]
\vskip3\jot
\centering
\noindent
\psfig{figure=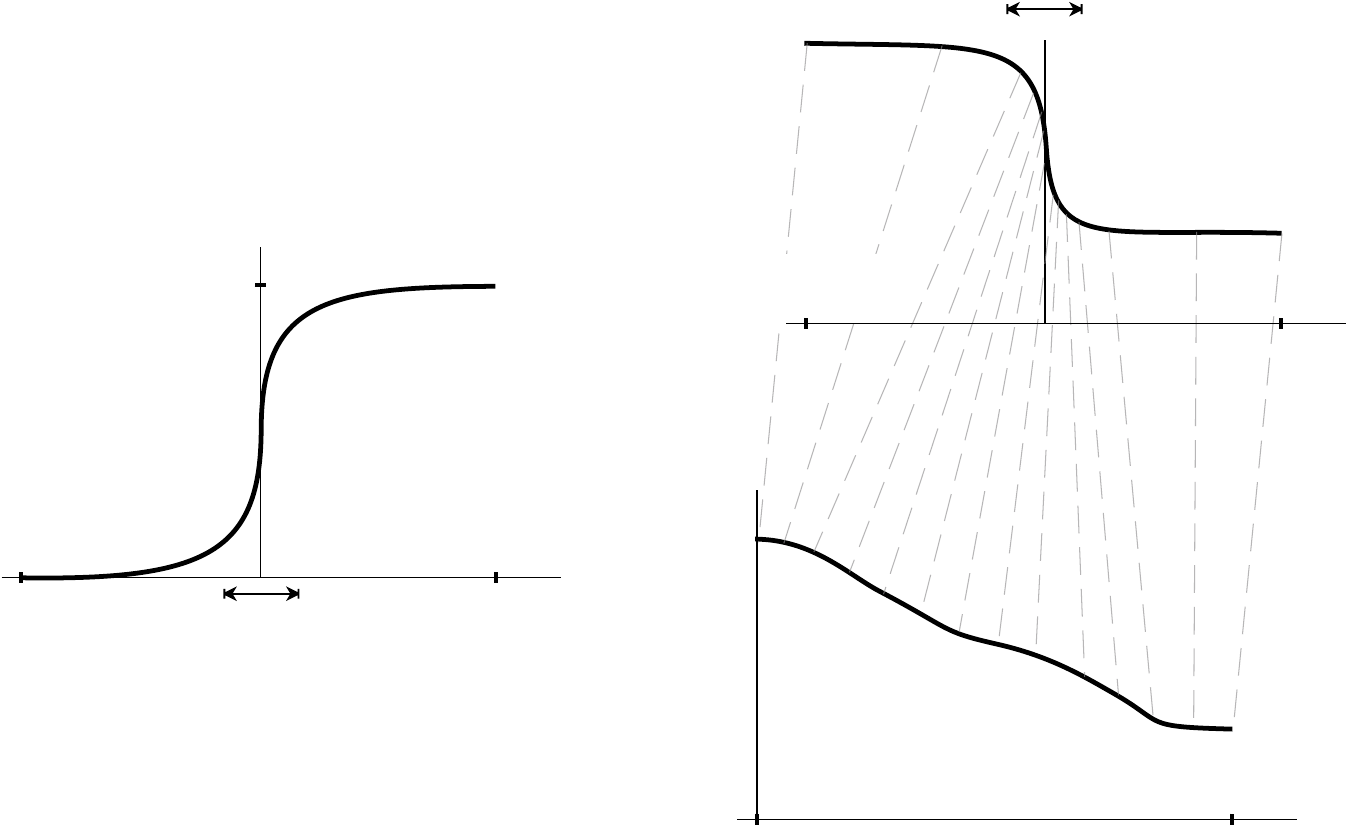,height=7cm}
\vskip3\jot
\caption{The transformations from $\xi$ to $s$ and from $u_\e(\xi)$ to $\hue(s)$. The left-hand graph shows the bijection between $\xi\in[-1,1]$ and $s\in[0,1]$. The right-hand graphs illustrate how the transformation~\eqref{def:hatu} expands the region around $\xi=0$ and converts the function $u_\e$ with a near-discontinuity around $\xi=0$ into a function $\hue$ that has a slope of order $O(1)$.}
\label{fig:transformation}
\end{figure}

\begin{remark}
\label{rem:desingularizing}
The desingularizing effect of the transformation from $u_\e$ to
$\hue$  can best be recognized in
the last term in~\eqref{def:hatJ}. In terms of $\hue$, this term is
the $H^1$-seminorm of $\surd\hue$, and indeed boundedness of
$\hat\J_\e$ implies boundedness of $\surd\hue$ in $L^2(0,T;H^1(\rmS))$
(see the proof of Theorem~\ref{th:compactness}). Compare this with the
corresponding term in the non-transformed version~\eqref{eq:45bis},
where the vanishing of $\tau_\e\gamma_\e$ close to $\xi=0$ allows for
large gradients at that point. 

As an independent way of viewing the effect of the transformation, the equation~\eqref{pb:main_u} for $u_\e$ transforms into the equation
\[
\hgdens\; \partial _t \hue =  \partial_{ss} \hue
\]
for $\hue$.
Here the structure of the term $\partial_s^2 \hue$ (specifically, the lack of singular terms \emph{inside} the derivatives) is related to the better behaviour of $\hue$.
\end{remark}

\subsection{Formal analysis}

The most complete understanding of the limit $\e\to0$, and how $\J_\e$ becomes converted into $\J_0$,  arises from considering the rescaled function $\hat u_\e$. However, one can also gain some understanding from a formal analysis in the original variable $\xi$, which we outline here.

Since $(\rho_\e,\rho_\e v_\e)\in \CE0T\Xi$, we have
\begin{equation}
\label{eq:CEu}
g_\e \partial_t u_\e + \partial_\xi (u_\e v_\e g_\e) = 0,
\end{equation}
which can be partially solved to find
\begin{equation}
\label{def:w_e}
w_\e(\xi) := u_\e(\xi) v_\e(\xi) g_\e(\xi)   = \int_{-1}^\xi \dot u_\e g_\e.
\end{equation}
Since $g_\e$ vanishes in the interior of $[-1,1]$, we expect that $w_\e$ becomes  constant in $\xi\in(-1,1)$, with value 
\begin{equation}
\label{def:w-formal}
\int_{-1}^0 \dot u_\e g_\e \approx \frac12 \dot u(-1,t) =: w.
\end{equation}
%We now show how these formal characterizations imply the convergence of $\J_\e$ to $J_0$. 
Writing the first term in~\eqref{def:JWass} in terms of $w_\e$, and using the convergence of $w_\e$ to the constant~$w$, we find that the two $\xi$-integrals become
\[
\frac1{2\tau_\e} \int_{-1}^1 \frac{w_\e^2}{u_\e g_\e} \, d\xi
+ \frac{\tau_\e}2 \int_{-1}^1 \frac{\partial_\xi u_\e^2}{u_\e} g_\e \, d\xi
\ \stackrel{\e\to0} \approx\  \frac{w^2}{2} \int_{-1}^1 \frac1{u_\e}\, \frac{d\xi}{\tau_\e g_\e}
+ \frac12\int_{-1}^1 \frac{\partial_\xi u_\e^2}{u_\e} \, \tau_\e g_\e \, d\xi.
\]
This expression is already close to the definition of the function $M$ in~\eqref{def:M}. The remaining difference, the two factors $\tau_\e g_\e$ in the two integrals, is removed by the transformation~\eqref{diffeo} from $\xi$ to $s$. 

One can  go further, and reconstruct how the limit equation~\eqref{pb:limit} arises in the $\xi$-variable. 
For small $\e$, the large value of $\tau_\e$ implies that for fixed $t$ the function $\xi\mapsto u_\e(\xi,t)$ is close to a stationary solution of~\eqref{pb:main_u} with Dirichlet boundary conditions at $\xi=\pm1$. Writing this equation as
\begin{equation}
\label{pb:main_u2}
g_\e \partial_t u_\e = \tau_\e \partial_\xi \bigl( g_\e \partial_\xi u_\e\bigr),
\end{equation}
it follows that (writing $u_\pm$ for the boundary values)
\begin{equation}
\label{approx:dxiu}
\partial_\xi u_\e(t,\xi) \approx \frac{u_+ - u_-}{2\kappa\tau_\e}  \frac1{g_\e(\xi)},
\end{equation}
where we have used the characterization~\eqref{def:tau_e} of $\tau_\e$. Comparing~\eqref{pb:main_u2} with~\eqref{eq:CEu} we remark that
\[
w_\e := u_\e v_\e g_\e  = \tau_\e g_\e \partial_\xi u_\e,
\]
so that 
\[
\frac12 \dot u_- \ \stackrel{\eqref{def:w-formal}}= \ w \ \approx\  w_\e\  =\  \tau_\e g_\e \partial_\xi u_\e
\ \stackrel{\eqref{approx:dxiu}}\approx \ \frac{u_+-u_-}{2\kappa} \ \stackrel{\eqref{eq:51}}= \ \frac12  \, k (u_+-u_-),
\]
which coincides with~\eqref{pb:limita}.

\section{Compactness}
\label{sec:compactness}

The main results of this section, Theorems~\ref{th:compactness-bounded} and~\ref{th:compactness},
describe compactness properties of sequences {$\rho_\e$}, and their transformed versions $\hat\rho_\e$,
for which the initial energy $\E_\e(\rho_\e(0))=\hE_\e(\hat\rho_\e(0))$ and {$\J_\e(\rho_\e)=\hJ_\e(\hat\rho_\e)$} are bounded.
%{The limit behaviour of $\rho_\e$ and of its rescaled version $\hat\rho_\e$}
%differ and get separate treatment.

Let us first comment on what one might expect. 
For {$\rho_\e$} and $\hat \rho_\e$ the limit objects are
measures $\rho_0$ and $\hat \rho_0$ concentrated in $\{-1,1\}$ and $\{-\kappa,\kappa\}$.
The existence of converging subsequences is a simple consequence of the bounded total variation of the measures and the bounded domain of definition. However, this convergence alone does not contain enough information for the lower bound result of Theorem~\ref{th:inner_lower_bound}. 

The key to obtaining more detailed convergence statements lies in using the objects that appear in $\J_\e$ and $\hJ_\e$, which are $w_\e := v_\e u_\e g_\e$, $\hat w_\e$, $\partial_\xi u_\e$, and $\partial_s \hue$.
Boundedness of $\J_\e(\rho_\e)$ implies, using the definitions~\eqref{def:JWass} and~\eqref{def:hatJ}, that there exists a constant $C>0$ such that
\begin{alignat}2
\frac1{\tau_\e}&\int_0^T\!\!\int_{\Xi} \frac{w_\e^2}{u_\e g_\e}\,\d \xi\d t \leq C,
&\qquad
\tau_\e&\int_0^T\!\!\int_\Xi \frac{\partial_\xi u_\e^2}{u_\e}\, \dmm{\gamma_\e} \xi \d t \leq C,
\label{bound:c0}\\
&\int_0^T\!\!\int_{\rmS} \frac{\hat w_\e^2}{\hat u_\e}\,\d s\d t \leq C,
&\qquad
&\int_0^T\!\!\int_\rmS \frac{\partial_s \hue^2}\hue\, \d s \d t \leq C.
\label{bound:c1}
\end{alignat}
If the sequence $\hue$ also happens to be bounded in $L^\infty$, then the bounds~\eqref{bound:c1} imply weak compactness of $\hat w_\e$ and $\partial_s\hue$ in $L^2((0,T)\times\rmS)$. Also, since $\tau_\e\gamma_\e$ is unbounded in any set not containing $\xi=0$, \eqref{bound:c0} suggests that $u_\e$ should become constant in $[-1,0)$ and $(0,1]$. In Theorem~\ref{th:compactness-bounded}, where we make this additional assumption of boundedness in $L^\infty$, we show that the remarks above indeed are true.
Moreover, we shall see that we can recover the canonical decomposition
$\rho_0=\frac 12 u^-\delta_{-1}+\frac 12u^+\delta_1$ by
taking the limit of the traces of the densities $u_\e$ at $\xi=\pm1$, and similarly for $\hat\rho_0$.

When $\hat u_\e$ is not assumed to be bounded in $L^\infty$, singular behaviour is possible that violates the $L^\infty$ bound but influences neither the boundedness of energy and dissipation nor the limit object~$\rho_0$. We treat this case in Theorem~\ref{th:compactness}.

\begin{theorem}[Compactness under uniform $L^\infty$ bounds]
\label{th:compactness-bounded}
Let $\rho_\e=u_\e\gamma_\e\in C([0,T];\Z^\text{meas})$
satisfy, for suitable constants $\mass,C>0$ and for all
$\e>0$, the estimates
\begin{equation}
  \label{bounds:th:compactness-bounded}
 \rho_\e(t,\Xi)=\mass 
  \text{ for all } t\in[0,T],\quad
  \E_\e(\rho_\e(0))\le C,\quad
  \J_\e(\rho_\e)\le C, \quad\text{and}\quad
  \|u_\e\|_\infty \leq C.
\end{equation}
Then there exists a subsequence (not relabeled) and a limit
{$\rho_0=u_0\gamma_0\in
C([0,T];\Z^\text{meas})$}
such that the following hold:
\begin{enumerate}
\item \label{th:compactness-b:MFconv}
   $\rho_\e(t) \longweakstarto \rho_0(t)$ in $\M(\Xi)$
  for every $t\in [0,T]$;
 
\item \label{th:compactness-b:traceconv}
  The spatial traces $u_\e(\cdot , \pm 1)$ are well-defined and
  converge {strongly} to $u_0^\pm(\cdot)= u_0(\cdot,\pm1)$ in
  {$L^1(0,T)$};
  
\item \label{th:compactness-b:locconstant}
For all $0<\delta<1$ the function $u_\e$ converges uniformly to the limiting trace values $u_0^\pm$ in $L^1(0,T;L^\infty(-1,-\delta))$ and $L^1(0,T;L^\infty(\delta,1))$. 
\end{enumerate}

\medskip

Let $\hue$ be the transformed sequence and let $\hat w_\e$ be given
as in \eqref{def:hatJ}.
Then there exist limits
$\huz\in L^\infty((0,T)\times\rmS)\cap L^2(0,T;W^{1,2}(\rmS))$ and $\hat w_0\in
L^2((0,T)\times \rmS)$ such that 
\begin{enumerate}[resume]
\item\label{th:compactness-b:L2H1}
$\hue\longweakstarto \huz$ in $L^\infty((0,T)\times\rmS)$, and $\partial_s\hue \weakto \partial_s \huz$,
$\hat w_\e\weakto \hat w_0$
in $L^2((0,T)\times \rmS)$;
\item\label{th:compactness-b:identificationoftraces}
  the traces $\hat u_0^\pm$ of $\hat u_0(\cdot,s)$ at $s=\pm \kappa$
  coincide with the traces of $u_0$, i.e.
  they satisfy $\hat u_0^\pm=u_0^\pm$ in $(0,T)$;
\item\label{th:compactness-b:hatMF}
  $\hat\rho_\e\longweakstarto 
\hat\rho_0(t)=\frac 12(
u_0^-(t)\delta_{-\kappa}+u_0^+(t)\delta_\kappa)$ in $\M(\rmS)$ for
every $t\in [0,T]$;
\item\label{th:compactness-b:continuityequation}
$\hat w_0(t,\cdot)$ is constant in $\rmS$ for a.e.\
$t\in (0,T)$ and satisfies $\hat w_0(t,\cdot)= \frac 12\dot u_0^+(t)$
a.e.\ in $(0,T)$.
The couple $(\hat\rho_0,\hat \nu_0)$, $\hat \nu_0=\hat w_0\L^2\restr{(0,T)\times\rmS}$ satisfies the
continuity
equation
\begin{equation}
  \label{th:compactness-b:continuityeq}
  \partial_t\hat\rho_0+\partial_s \hat\nu_0=0\quad\text{in }\mathscr D'((0,T)\times\R).
\end{equation}
\end{enumerate}
\end{theorem}

We will see in Theorem~\ref{th:compactness_solutions} that in the
special case of \emph{solutions}
of~\eqref{pb:main}, which satisfy $\hat\A_\e(\hat \rho_\e)=0$
rather than $\hat \A_\e(\hat\rho_\e)\leq C$,
the limit object $\huz$ is a linear interpolation of the values at $s=\pm\kappa$.

\begin{proof}

We divide the proof in a few steps; we will denote by $C$ various
constants which are independent of $\e$.

\smallskip

\textbf{Step 1: Entropy estimates.} \emph{There exists a constant
  $C>0$ such that}
\begin{equation}
  \label{bound:energy}
  \E_\e(\rho_\e(t))=\hat \E_\e(\hat
  \rho_\e(t))\le C\quad\text{\emph{for every }}t\in [0,T];
\end{equation}
\emph{in particular, for any  subsets $A\Subset (-1,1)$ and $\hat
  A\Subset (-\kappa,\kappa)$ we have}
\begin{equation}
  \label{vanishing-rho}
  \lim_{\e\to0}\sup_t \rho_\e(t,A)=0,\quad
  \lim_{\e\to0}\sup_t \hat\rho_\e(t,\hat A)=0.
\end{equation}
It is sufficient to prove \eqref{bound:energy} for the unrescaled measures
$\rho(t)$.
First we note that $\E_\e$ is nonnegative;
denoting by $\check\rho(t):=\rho(T-t)$ the time-reversed curve,
since $\J_\e(\check\rho)=\J_\e(\rho)$,
the bounds \eqref{bounds:th:compactness-bounded} and \eqref{cond:nonnegativeA}
imply that $\E_\e(\rho_\e(t))\le C$ for every $t\in
[0,T]$.

Property~\eqref{vanishing-rho} follows from \eqref{bound:energy} and the fact that
$\lim_{\e\downarrow0}\gamma_\e(A)=\lim_{\e\downarrow0}\hat\gamma_\e(\hat
A)=0$.
Considering e.g.\ the case of $A\Subset(-1,1)$,
by using first the inequality $r\log r\ge -\rme^{-1}$ and then Jensen
inequality, we get for every $A\Subset (-1,1)$
\begin{equation}
  \label{eq:50bis}
  \frac 1\rme\gamma_\e((-1,1)\setminus
  A)+\E_\e(\rho_\e(t))+\mass\log\mass
  \ge \int_A u_\e(t,\xi)\log u_\e(t,\xi)\,\dmm{\gamma_\e}\xi\ge \rho_\e(t,A)\log\Big(\frac{\rho_\e(t,A)}{\gamma_\e(A)}\Big),
\end{equation}
so that $\gamma_\e(A)\to 0$ implies $\sup_t \rho_\e(t,A)\to0$ as
$\e\to0$.

\smallskip

\textbf{Step 2: Estimates on \boldmath $\hat u_\e, \hat w_\e$.}
\emph{There exists a constant $C>0$ such that }
\begin{equation}
  \label{bound:winL2}
  \|\hue\|_\infty \leq C\quad \text{and}\quad 
  \int_0^T \!\!\int_\rmS \hat w_\e(t,s)^2\,\d s\d t\le C.
\end{equation}
The first bound derives from assumption
\eqref{bounds:th:compactness-bounded}, and the second follows
easily from~\eqref{bound:c1} and the $L^\infty$-bound on $\hat
u_\e$. We state these here to contrast with the corresponding, weaker,
versions in the proof of Theorem~\ref{th:compactness}.

\smallskip

\textbf{\boldmath Step 3: Pointwise weak convergence of $\hat \rho_\e(t)$ (statement \ref{th:compactness-b:hatMF})}:
\emph{there exists a sequence $\e_n\downarrow0$ and a
  limit $\hat \rho_0(t)\ll \hgz$ such
  that $\hat\rho_{\e_n}(t)\weakstarto \hat \rho_0(t)$ for every $t\in
  [0,T]$.} 

Starting from the continuity equation
we have for every $0\le
t_0<t_1\le T$ and $\varphi\in C^1([0,1])$,
%\begin{equation}
%  \label{eq:60}
\[
  \int_\rmS \varphi\,\d\hat\rho_{\e}(t_1)-
  \int_\rmS \varphi\,\d\hat\rho_{\e}(t_0)=
  \int_{t_0}^{t_1}\!\! \int_\rmS\partial_s \varphi\; \hat w_\e\,\d s\d t.
\]
Recalling
the definition of the $L^1$-Wasserstein distance $\sfd_{W_1}$~\cite{AmbrosioGigliSavare05}, we find
\begin{align*}
  \sfd_{W_1}(\rho_{\e}(t_1),\rho_{\e}(t_0))&:=
  \sup\biggl{\{}  \int_\rmS \varphi\,\d\hat\rho_{\e}(t_1)-
  \int_\rmS \varphi\,\d\hat\rho_{\e}(t_0):\varphi\in C^1([0,1]),\
  |\partial_s \varphi|\le 1\biggr{\}}\\&\le
  \int_{t_0}^{t_1}\!\!\int_\rmS | \hat w_\e(t,s)|\,\d s\d t
  \leq \sqrt {2\kappa (t_1-t_0)} \biggl{(}\int_0^T \!\!\int_\rmS \hat w_\e(t,s)^2\, \d s\d t\biggr{)}^{1/2}.
\end{align*}
It follows by~\eqref{bound:winL2} that the curves $t\mapsto \hat\rho_{\e,t}$ are an equicontinuous family of mappings from $[0,T]$ into the space $\M(\rmS)$ endowed with
the $L^1$-Wasserstein distance; since the total mass is $\mass$, the claim follows by
the Arzel\`a-Ascoli
theorem.

\smallskip

\textbf{\boldmath Step 4. Weak convergence of $\rho_\e(t)$ (statement \ref{th:compactness-b:MFconv}).}
\emph{Writing the limit $\hat \rho_0$ of the previous step as $\hat \rho_0(t)=\frac 12 u_0^-(t)\delta_{-\kappa}+\frac
12u_0^+(t)\delta_\kappa$, we have for every $t\in[0,T]$,}
\begin{equation}
  \label{conv:rho}
  \rho_{\e_n}(t)\weakstarto \rho_0(t)=\frac 12 u_0^-(t)\delta_{-1}+\frac 12
  u_0^+\delta_1.
\end{equation}
Let us fix $t\in [0,T]$ and let us consider a subsequence $\e_n'$ of
$\e_n$ along which $\rho_{\e_n'}(t)$ converges weakly to some
$\tilde \rho\in \M(\Xi).$
By the result of Step 1, we know that $\tilde \rho=\frac 12
\tilde u^-\delta_{-1}+\frac 12 \tilde u^+\delta_1$ for some $\tilde u^\pm$: if we show that
$\tilde u^\pm=u^\pm_0(t)$ we have proved the thesis.
Considering $\tilde u^-$ and taking a function $\hat\varphi\in
C(\rmS)$, we know that $\varphi_\e:=\hat\varphi\circ s_\e$
converges pointwise to $\varphi_0(\xi)=\hat\varphi\circ s_0$, where
$s_0(\xi)=\mathrm{sign}(\xi)$ (with $s_0(0)=0$)
and the convergence is uniform on compact subsets of $[-1,1]\setminus
\{0\}$.
It then follows that
\begin{align*}
  \frac 12\sum_\pm \hat u^\pm_0(t)\hat\varphi(\pm\kappa)&=
  \lim_{n\to\infty}\int_\rmS \hat\varphi(s)\,\hat\rho_{\e_n'}(t,\d s)=
  \lim_{n\to\infty}\int_1^1 \varphi_{\e_n'}(\xi)\,\rho_{\e_n'}(t,\d\xi)=
  \int_{-1}^1 \varphi_0(\xi)\,\tilde\rho(\d\xi)
  \\&=
  \frac 12\sum_\pm\tilde u^\pm\hat\varphi(\pm\kappa).
\end{align*}
Since $\hat\varphi$ is arbitrary, we obtain~\eqref{conv:rho}.

This shows the intuitive result that the densities of $u_0$ and $\hat u_0$ (with respect to $\gamma_0$ and $\hgz$) are the same; we call them both $u_0^\pm$. It mirrors the fact that the traces of $u_\e$ (in $\xi = \pm1$) and of $\hat u_\e$ (in $s=\pm \kappa$) are also the same. 

\smallskip

\textbf{Step 5. Convergence of \boldmath $u_\e$ (statements
  \ref{th:compactness-b:traceconv} and~\ref{th:compactness-b:locconstant}):}
\emph{The traces $u_{\e_n}^\pm=u_{\e_n}(\,\cdot\,,\pm 1)$ strongly converge in $L^1(0,T)$ to the limits $u_0^\pm$ defined in the previous step. In addition, setting $\omega^\pm_\delta:= \pm (\delta,1)$ we have}
\begin{equation}
  \label{eq:66}
  \lim_{n\to\infty}\int_0^T  \sup_{\xi\in \omega^\pm_\delta}|u_{\e_n}(t,\xi)-
   u_0^\pm(t)|\, \d t = 0 \quad
   \text{for every }0<\delta<1.
\end{equation}
Let us first observe that the quantities 
\begin{displaymath}
  \bar
  u^\pm_{\e_n,\delta}(t):=\frac{\rho_{\e_n,t}(\omega_\delta^\pm)}{\gamma_\e (\omega_\delta^\pm)}=
  \frac 1{\gamma_\e
    (\omega_\delta^\pm)}\int_{\omega_\delta^\pm}u_{\e_n}(t,\xi)\,\d\gamma_\e(\xi),  
  \quad 0<\delta<1,
\end{displaymath}
are uniformly bounded and converge pointwise to $u^\pm_0(t)$ for every
$t\in (0,T)$ by Step 4.
By~\eqref{bound:c0} and the boundedness of $u_\e$ we have
\[
  \label{eq:52bis}
  \lim_{\e\downarrow0}\int_0^T |\theta_\e^\pm|^2(t)\,\d t=0,\quad\text{where}\quad
  \theta_\e^\pm(t):=\sup_{\xi,\eta\in
    \omega^\pm_\delta}|u_\e(t,\xi)-u_\e(t,\eta)|\le
  \Big(\delta\int_{\omega^\pm_\delta}|\partial_\xi u_\e|^2\,\d\xi\Big)^{1/2}.
\]
We then calculate for $\xi,\eta\in\omega_\delta^\pm$, 
\[
|u_\e(t,\xi) - \bar u^\pm_{\e_n,\delta}(t)|
\leq \frac1{\gamma_\e(\omega_\delta^\pm)} \int_{\omega_\delta^\pm} 
  |u_\e(t,\xi) - u_\e(t,\eta)|\,\dmm{\gamma_\e}\eta
\leq \frac{\theta_\e^\pm(t)}{\gamma_\e(\omega_\delta^\pm)}.
\]
Recalling that $\gamma_\e(\omega_\delta^\pm)\to 1/2$ as
$\e\to0$, we thus obtain
\[
  \label{eq:65}
  \lim_{n\to\infty}\int_0^T     \sup_{\xi\in \omega^\pm_\delta}|u_{\e_n}(\xi,t)-
  \bar u_{{\e_n},\delta}^\pm(t)|\,\d t=0
\]
%\end{equation}
which in particular yields \eqref{eq:66} and the convergence of the traces, since
$\bar u_{{\e_n},\delta}^\pm$ strongly converge to $u_0^\pm$ in every
$L^p(0,T)$, $p<+\infty$.

\smallskip

\textbf{Step 6. Compactness and limits (statement~\ref{th:compactness-b:L2H1}).}
Given the estimates~\eqref{bound:winL2} and~\eqref{bound:c1}, this follows from standard results. 

\smallskip

\textbf{Step 7. Identification of the traces of \boldmath $\huz$  (statement \ref{th:identificationoftraces}).}
Since the trace operator $\Tr$ is weakly continuous from $H^1(\rmS)$ to $L^2(\{-\kappa,\kappa\})\simeq \R^2$, the weak convergence of $\hat u_{\e_n}$ in $H^1(\rmS)$ implies that its traces $\hat u_{\e_n}(\cdot,\pm\kappa) = u_{\e_n}(\cdot,\pm 1)$ converge weakly in $L^2(0,T;\R^2)$ to $\Tr \hat u_0$. Since $u_{\e_n}(\cdot,\pm 1)$ converges strongly to $u_0^\pm$ in $L^1(0,T)$ (Step 5), it follows that $\Tr \hat u_0 = u_0^\pm$.

\smallskip
\textbf{Step 8: The continuity equation and the structure of \boldmath $\hat w_0$ (statement \ref{th:compactness-b:continuityequation})}.
Passing to the limit in the continuity equation
\eqref{eq:42bis} and using the previous convergence result we immediately
find \eqref{th:compactness-b:continuityeq}.
Since $\hat \rho_0$ is supported in $[0,T]\times \{-\kappa,\kappa\}$,
we obtain that $\partial_s \hat \nu_0=0$ in $[0,T]\times (-\kappa,\kappa)$, so that $w_0$ depends only on $t$. 

Choosing a test function of the form $\varphi(t,s)=\psi(t) \zeta(s)$
with $\psi\in C^\infty_{\rm c}(0,T)$ and $\zeta\in C^\infty_{\rm
  c}(\R)$, $\zeta\equiv 1$ on a neighborhood of $\kappa$ and
$\zeta\equiv 0$ on $(-\infty,0]$, we
obtain
\[
    \frac 12 \int_0^T \dot\psi(t)\hat u_0^+(t)\,\d t
    =\int_0^T\!\!\int_\R \dot\psi(t)\zeta(s) \,\hat\rho_0(t,\d s)\,\d t=
    -\int_0^T\!\!\int_\rmS \psi(t)\zeta'(s)\hat w_0(t)\,\d s\,\d t.
\]
Since $\int_\rmS
  \zeta'(s)\,\d s=\zeta(\kappa)=1$, we conclude that $\hat w_0$ is the
  distributional derivative of
  $\frac 12\hat u_0^+$.
This concludes the proof of Theorem~\ref{th:compactness-bounded}.
\end{proof}

\medskip

We now discuss the case where $\hue$ is not assumed to be bounded in $L^\infty$. A simple example shows how $\hue$ may become singular without affecting any of the relevant limit processes. Take any sequence $\hat \rho_\e$ with bounded $\hE_\e(\hat\rho_\e(0))$ and $\hJ_\e(\hat\rho_\e)$, and let $\hat u_\e=\d\hat\rho_\e/\d\hat\gamma_\e$ be bounded from above and away from zero as well. Fix two nonnegative functions $\varphi\in C^1_c((-\kappa,\kappa))$  and $\psi\in C^1_c(\R)$, fix $0<t_0<T$, and define
\[
\widetilde\rho_\e(t,s) := \hat\rho_\e(t,s) + \frac1{\sqrt\e} \psi\Bigl(\frac{t-t_0}\e\Bigr)\varphi(s)\hat\gamma_\e(s).
\]
Note that since the additional term blows up polynomially, while $\hat \gamma_\e$ converges to zero exponentially fast on $\supp \varphi$, the limits of $\widetilde\rho_\e$ and $\hat \rho_\e$ are the same. For the same reason the perturbed $\widetilde w_\e$, satisfying $\partial_t\widetilde\rho_\e + \partial_s\widetilde w_\e=0$, only differs from $\hat w_\e$ by an exponentially small amount. Therefore 
\[
\limsup_{\e\to0} 
\int_0^T \!\! \int_\rmS \frac{\widetilde w_\e^2}{\widetilde u_\e}\, \d s\d t
= \limsup_{\e\to0} 
\int_0^T \!\! \int_\rmS \frac{\hat w_\e^2}{\hat u_\e}\, \d s\d t
<\infty.
\]
We also estimate
\[
\int_0^T \!\!\int_\rmS \frac{\partial_s \widetilde u_\e^2}{\widetilde u_\e}
\, \d s\d t
\leq 2\int_0^T \!\!\int_\rmS \frac{\partial_s \hat u_\e^2}{\hat u_\e}
\, \d s\d t
+ \frac 2\e \int_0^T \psi^2\Big(\frac{t-t_0}{\e}\Bigr) 
    \int_\rmS \frac{\varphi'^2}{\hat u_\e}\, \d s\d t.
\]
The first term of this estimate is bounded  by assumption, and the second is bounded by the scaling in $\e$ and the assumption that $\hat u_\e$ is uniformly bounded away from zero. 

This example shows that the assumptions of bounded initial energy and bounded $\J$ do not rule out singular behaviour of the sequence $\hat u_\e$ \emph{between} $-\kappa$ and $\kappa$. The example also suggests what form this singular behaviour might take: that of a singular measure in time (called $\lambda^\perp$ below), but with bounded total variation. This is exactly what we prove in the following theorem.

\begin{theorem}[Compactness, the general case]
\label{th:compactness}
Let $\rho_\e=u_\e\gamma_\e\in C([0,T];\Z^\text{meas})$
satisfy, for suitable constants $\mass,C>0$ and for all
$\e>0$, the estimates
\begin{equation*}
%  \label{eq:48bis}
  \rho_\e(t,\Xi)=\mass \text{ for all }t\in[0,T],
  \quad \E_\e(\rho_\e(a))\le C,\quad\text{and }
  \J_\e(\rho_\e)\le C.
\end{equation*}
Then there exists a subsequence (not relabeled) and a limit
{$\rho_0=u_0\gamma_0\in
C([0,T];\Z^\text{meas})$}
such that the following hold:
\begin{enumerate}
\item \label{th:compactness:MFconv}
   $\rho_\e(t) \stackrel*\longrightharpoonup \rho_0(t)$ in $\M(\Xi)$
  for every $t\in [0,T]$;
 
\item \label{th:compactness:traceconv}
  The spatial traces $u_\e(\cdot , \pm 1)$ are well-defined and
  converge {strongly} to $u_0^\pm(\cdot)= u_0(\cdot,\pm1)$ in
  {$L^1(0,T)$};
\item \label{th:compactness:locconstant} For all $0<\delta<1$ the
  function $u_\e$ converges uniformly to the limiting trace
  values $u_0^\pm$ in $L^1(0,T;L^\infty(-1,-\delta))$ and
  $L^1(0,T;L^\infty(\delta,1))$.
\end{enumerate}

\medskip

Let $\hue$ be the transformed sequence and let $\hat w_\e$ be given
as in \eqref{def:hatJ}.
Then there exist limit functions
$\huz\in L^1(0,T;W^{1,1}(-\kappa,\kappa))$, $\hat w_0\in
L^1((0,T)\times\rmS)$, a reference singular measure $\lambda^\perp\in
\M([0,T])$ with $\lambda^\perp \perp \L^1$, and a function
$\hat m_0\in L^\infty_{\Lambda^\perp}([0,T]\times \rmS)$ with $\partial_s \hat m_0 \in L^2_{\Lambda^\perp}([0,T]\times \rmS)$, where
$\Lambda^\perp=\lambda^\perp\otimes \L^1\restr{\rmS} \in \M([0,T]\times
\rmS)$, such that 
\begin{enumerate}[resume]
\item \label{th:compactness:L2H1alex}
$\hue\longweakstarto \huz+\hat m_0 \Lambda^\perp$,
$\partial_s\hue \longweakstarto \partial_s
\huz+\partial_s \hat m_0 \Lambda^\perp$, and
$\hat w_\e\longweakstarto \hat w_0$
in the duality with $C([0,T]\times\rmS)$;
\item\label{th:identificationoftraces}
  the traces $\hat u_0^\pm$ of $\hat u_0(\cdot,s)$ at $s=\pm \kappa$
  coincide with the traces of $u_0$, i.e.
  they satisfy $\hat u_0^\pm=u_0^\pm$ a.e. in $(0,T)$;
  the traces $\hat m_0^\pm(t)$ of $\hat m_0(t,\cdot)$ vanish
    for $\lambda^\perp$-a.e.\ $t\in [0,T]$;
\item\label{th:compactness:hatMF}
  $\hat\rho_\e\stackrel*\longrightharpoonup 
\hat\rho_0(t)=\frac 12(
u_0^-(t)\delta_{-\kappa}+u_0^+(t)\delta_\kappa)$ in $\M(\rmS)$ for
every $t\in [0,T]$;
\item\label{th:continuityequation}
$\hat w_0(t,\cdot)$ is constant in $\rmS$ for a.e.\
$t\in (0,T)$ and satisfies $\hat w_0(t,\cdot)= \frac 12\dot u^+(t)$
a.e.\ in $(0,T)$.
The couple $(\hat\rho_0,\hat \nu_0)$, $\hat \nu_0=\hat w_0\L^2\restr{(0,T)\times\rmS}$ satisfies the
continuity
equation
\begin{equation}
  \label{eq:75}
  \partial_t\hat\rho_0+\partial_s \hat\nu_0=0\quad\text{in }\mathscr D'((0,T)\times\R).
\end{equation}
\end{enumerate}
\end{theorem}

\begin{remark}
Parts~\ref{th:compactness:MFconv}--\ref{th:compactness:locconstant}, \ref{th:compactness:hatMF}, and~\ref{th:continuityequation} are the same as in Theorem~\ref{th:compactness-bounded}. The main difference lies in the structure of the limits of $\hue$ and $\hat w_\e$ (part~\ref{th:compactness:L2H1alex}) and therefore the identification of the traces of $\huz$ and $\hat m_0$ (part~\ref{th:identificationoftraces}).
\end{remark}

\begin{proof}
Some of the steps are the same as in the case of Theorem~\ref{th:compactness-bounded}; for those we only give the statement. For the others we detail the differences. 

\smallskip

\textbf{Step 1: Entropy estimates.} \emph{There exists a constant
  $C>0$ such that}
\begin{equation*}
  \label{eq:55}
  \E_\e(\rho_\e(t))=\hat \E_\e(\hat
  \rho_\e(t))\le C\quad\text{\emph{for every }}t\in [0,T];
\end{equation*}
\emph{in particular, for any  subsets $A\Subset (-1,1)$ and $\hat
  A\Subset (-\kappa,\kappa)$ we have}
\begin{equation*}
  \label{eq:56}
  \lim_{\e\to0}\sup_t \rho_\e(t,A)=0,\quad
  \lim_{\e\to0}\sup_t \hat\rho_\e(t,\hat A)=0.
\end{equation*}

\smallskip

\textbf{Step 2: Estimates on \boldmath $\hat u_\e,\hat w_\e$.}
\emph{There exists a constant $C>0$ such that }
\begin{equation}
  \label{eq:57}
  \int_0^T \sup_{s\in \rmS}|\hat u_\e(t,s)|\,\d t\le C,\quad
  \int_0^T \!\!\int_\rmS |\hat w_\e(t,s)|\,\d s\,\d t\le C.
\end{equation}
These bounds are weaker than the $L^\infty$-bound on $\hat u_\e$ and the $L^2$-bound on $\hat w_\e$ of Theorem~\ref{th:compactness-bounded}.
Let us set $\hat p_\e:=\sqrt {\hat u_\e}$: since $\hat p_\e\in
L^2(0,T;W^{1,2}(-\kappa,\kappa))$ its traces at $s=\pm\kappa$ are well defined
and belong to $L^2(0,T)$.
We set
$\hat \theta_\e(t):=\sup_{r,s\in \rmS}|\hat
p_\e(t,r)-\hat p_\e(t,s)| $.
Standard estimates yield
\[
%  \label{eq:58}
  \hat \theta_\e(t)\le
  \Big(2\kappa\int_\rmS|\partial_s \hat p_\e|^2\,\d s\Big)^{1/2},\quad
  \int_0^T \hat\theta_\e^2(t)\,\d t\le \hat\J_\e(\hat\rho_\e;0,T)\le C.
\]
Moreover
\[
%  \label{eq:59}
  \Big(\int_\rmS \hat p_\e\,\d\hat\gamma_\e\Big)^2\le
  \int_\rmS \hat p_\e^2\,\d\hat\gamma_\e=\mass,
\]
and
%\label{eq:54bis}
\[
    \hat p_\e(t,s)\le
    \int_{\rmS}\hat p_\e(r)\,\d\hat \gamma_\e(r)+
    \hat\theta_\e(t)\le  \sqrt \mass+\hat\theta_\e(t)
    \quad \text{for every }s\in \rmS,
\]
and therefore
\[
    \sup_{s\in \rmS}\hat u_\e(t,s)\le
    2 \mass+2\hat\theta_\e^2(t).
\]
The second estimate of~\eqref{eq:57} then follows from
\begin{displaymath}
  \int_\rmS |\hat w_\e(t,s)|\,\d s\le \Big(\int_\rmS \frac{|\hat
    w_\e(t,s)|^2}{\hat u_\e}\,\d s\Big)^{1/2}\Big(\int_\rmS \hat
  u_\e(t,s)\,\d s\Big)^{1/2}.
\end{displaymath}

\smallskip

\textbf{\boldmath Step 3: Pointwise weak convergence of $\hat \rho_\e(t)$ (statement \ref{th:compactness:hatMF})}:
\emph{there exists a sequence $\e_n\downarrow0$ and a
  limit $\hat \rho_0(t)$ such
  that $\hat\rho_{\e_n}(t)\weakstarto \hat \rho_0(t)$ for every $t\in
  [0,T]$, and $\hat\rho_0(t)\ll \hgz$.}

\smallskip

\textbf{\boldmath Step 4. Weak convergence of $\rho_\e(t)$ (statement \ref{th:compactness:MFconv}).}
\emph{Writing the limit $\hat \rho_0$ of the previous step as $\hat \rho_0(t)=\frac 12 u_0^-(t)\delta_{-\kappa}+\frac12 u_0^+(t)\delta_\kappa$, we have for every $t\in[0,T]$,}
%\begin{equation}
\[
  %\label{eq:64}
  \rho_{\e_n}(t)\weakstarto \rho_0(t)=\frac 12 u_0^-(t)\delta_{-1}+\frac 12
  u_0^+\delta_1.
\]
%\end{equation}

\smallskip

\textbf{Step 5. Strong convergence of traces (statements
  \ref{th:compactness:traceconv} and~\ref{th:compactness:locconstant}):}
\emph{the traces
$u_{\e_n}^\pm=u_{\e_n}(\,\cdot\,,\pm 1)$ strongly converge in
$L^1(0,T)$ to the limits $u_0^\pm$ defined in the previous step. In addition, setting $\omega^\pm_\delta:= \pm (\delta,1)$ we have}
\begin{equation*}
  %\label{eq:66}
  \lim_{n\to\infty}\int_0^T  \sup_{\xi\in \omega^\pm_\delta}|u_{\e_n}(t,\xi)-
   u_0^\pm(t)|\, \d t = 0 \quad
   \text{for every }0<\delta<1.
\end{equation*}
The proof of this step is similar to that of Theorem~\ref{th:compactness}, but uses instead the estimate on $p_\e:=\sqrt {u_\e}$,
\[
%  \label{eq:52bis}
  \lim_{\e\downarrow0}\int_0^T |\theta_\e^\pm|^2(t)\,\d t=0,\quad
  \theta_\e^\pm(t):=\sup_{\xi,\eta\in
    \omega^\pm_\delta}|p_\e(\xi,t)-p_\e(\eta,t)|\le
  \Big(\delta\int_{\omega^\pm_\delta}|\partial_\xi p_\e|^2\,\d\xi\Big)^{1/2}.
\]

\textbf{Step 6. Compactness and limits (statement~\ref{th:compactness:L2H1alex}).} Because of the lack of an $L^\infty$ bound, from here on the proof differs significantly from that of Theorem~\ref{th:compactness-bounded}.
Let us set $\ell_\e(t):=1+\sup_{s\in \rmS}|\hat u_\e(t,s)|$.
Up to extracting a suitable subsequence $\e\to0$ (without changing notation)
we can
assume that there exist weak limits $\lambda\in \M([0,T])$ and $
\hat\mu_0,\hat\nu_0,\hat \varsigma_0,\hat\sigma_0\in \M([0,T]\times \rmS)$, such that (identifying functions with the corresponding measures)
\begin{equation}
  \label{eq:1}
  \ell_\e\longweakstarto \lambda,\quad
  \hat u_\e\longweakstarto \hat \mu_0,\quad
  \hat w_\e\longweakstarto \hat \nu_0,\quad
  \partial_s \hat u_\e\longweakstarto \hat \varsigma_0,\quad\text{and}\quad
  |\partial_s \hat u_\e|\longweakstarto \hat \sigma_0.
\end{equation}
Since $\hat u_\e(t,s)\le \ell_\e(t)$ we have $\hat \mu_0\le\Lambda := \lambda\otimes \L^1\restr{\rmS} $, so that $\hat \mu_0=\hat m_0 \Lambda$ for
a suitable bounded Borel
function $\hat m_0\in L_\Lambda^\infty([0,T]\times \rmS)$. Since
\begin{displaymath}
  \int_0^T \!\!\int_\rmS \Big(\frac{\hat w_\e^2}{\hat u_\e}+
  \frac{\partial_s \hat u_\e^2}{\hat u_\e}\Big)\,\d s\,\d t\le C,
\end{displaymath}
it follows (see Lemma~\ref{prop:L2bdd} below) that in
the limit $\hat \nu_0\ll \hat\mu_0$ and 
$|\hat \varsigma_0|\le \hat\sigma_0\ll\hat \mu_0$.  In particular
$\hat \nu_0=\hat n_0\Lambda$ and
$\hat\varsigma_0=\hat g_0 \,\Lambda$ with
$\hat n_0, \hat g_0 \in L_\Lambda^\infty([0,T]\times \rmS)$.
Since $\partial_s \hat \mu_0=\hat \varsigma_0$, we easily have for
every couple of test functions $\psi\in C^\infty([0,T]),\varphi\in
C^\infty_{\rm c}(-\kappa,\kappa)$,
\begin{displaymath}
  \iint_{(0,T)\times\rmS} \psi(t)\varphi(s)\,\hat g_0(t,s)\,\d
  s\,\dmm\lambda t=
  -\iint_{(0,T)\times\rmS} \psi(t)\varphi'(s)\,\hat m_0(t,s)\,\d
  s\,\dmm\lambda t.
\end{displaymath}
Since $\psi$ is arbitrary, we deduce that $\partial_s \hat m_0(t,\cdot)=
\hat g_0(t,\cdot)$ in $L^\infty(\rmS)$ for $\lambda$-a.e.\ $t\in [0,T]$.
We also deduce that
\begin{equation}
\label{eq:3}  \int_0^T \!\!\int_\rmS \Big(\frac{\hat n_0^2}{\hat m_0}+
  \frac{\hat g_0^2}{\hat m_0}\Big)\,\d s\,\dmm \lambda t\le C.
\end{equation}
The measure $\lambda \in \M([0,T])$ can be decomposed as 
 $\lambda=\ell\L^1+\lambda^\perp$ with $\ell \in L^1(0,T)$ and
 $\lambda^\perp \perp \L^1$, and similarly $\Lambda = \ell \L^1 \otimes \L^1\restr{\rmS} + \Lambda^\perp$ with $\Lambda^\perp = \lambda^\perp\otimes\L^1\restr{\rmS}$. We set
$\hat u_0:=\hat m_0\ell$ and
$\hat w_0:=\hat n_0\,\ell$, so that the limits in~\eqref{eq:1} can be decomposed as $\hat \mu_0 = \huz + \hat m_0\Lambda^\perp$, $\varsigma_0 = \partial_s \huz + \partial_s\huz \Lambda^\perp$, and $\hat \nu_0 = \hat w_0 + \hat n_0\Lambda^\perp$. In Step~8 below we show that the last term, $\hat n_0\Lambda^\perp$, vanishes.

\medskip

\textbf{Step 7. \boldmath $\hat \mu_0$ and
  $\hat u_0$ have equal traces
  (statement \ref{th:identificationoftraces}).}
Let us consider, e.g., the case of $-\kappa$ and
take nonnegative test functions $\psi\in C([0,T])$ with $\sup_{t\in
  [0,T]}|\psi(t)|\le 1$,
and $\varphi\in C([-\kappa,\kappa])$ with support in
$[-\kappa,0)$ and integral $1$, so that
$\Phi(s):=\int_s^\kappa\varphi(r)\,\d r$ is decreasing, supported in
$[-\kappa,0)$,
and satisfies $\Phi(-\kappa)=1$.
We also set
$\varphi_\delta(s):=\delta^{-1}\varphi(-\kappa+\delta^{-1}(s+\kappa))$,
$\Phi_\delta(s):=\Phi(-\kappa+\delta^{-1}(s+\kappa))$, $\varphi_\delta(s)=-\Phi_\delta'(s)$.
Denoting by $\Omega$ the product $[0,T]\times\rmS$, we have
\begin{align*}
  \Big|&\int_0^T \psi(t) u^-_0(t)\,\d t-
  \iint_{\Omega} \psi(t)\varphi_\delta(s) \,\hat \mu_0(\d s\d t)\Big|\le
  \int_0^T \big| u^-_0(t)- \hat u^-_\e(t,-\kappa)\big|\,\d t\\&+
  \iint_{\Omega} \psi(t)\varphi_\delta(s)|\hat u^-_\e(t,-\kappa)-\hat u_\e(t,s)|\,\d s\,\d t+
  \Big|\iint_{\Omega} \psi(t)\varphi_\delta(s)\hat
  u_\e(t,s)\,\d s\,\d t-
  \iint_{\Omega} \psi(t)\varphi_\delta(s)\,\hat \mu_0(\d s\d t)
  \Big|.
\end{align*}
Passing to the limit as $\e\to0$ the first and third terms vanish;
concerning the
second term, we have
\begin{align*}
  \iint_{\Omega} \psi(t)\varphi_\delta(s)|\hat u_\e(t,-\kappa)-
  \hat u_\e(t,s)|\,\d
  s\,\d t&\le
  \iint_{\Omega} \psi(t)\varphi_\delta(s)
  \Big(\int_{-\kappa}^s|\partial_s
  \hat u_\e(t,r)|\,\d r\Big)
  \,\d s\,\d t
  \\&=
  \int_0^T \psi(t)\Big(\int_{-\kappa}^\kappa \Phi_\delta(s) |\partial_s \hat
  u_\e(t,s)|\,\d s\Big)\,\d t .
\end{align*}
Combining these inequalities and passing to the limit as $\e\to0$
we get
\begin{align*}
  \Big|&\int_0^T \psi(t) u^-_0(t)\,\d t-
  \iint_{\Omega} \psi(t)\varphi_\delta(s) \,\hat \mu_0(\d s\d t)\Big|\le
  \iint_\Omega \psi(t)\Phi_\delta(s)\,\hat\sigma_0(\d s\d t),
\end{align*}
so that, applying Lebesgue's dominated convergence theorem with 
the fact that $\Phi_\delta(s)\to 0$ for $s>-\kappa$ and $0\leq\Phi_\delta\leq 1$, we obtain 
\begin{align*}
  \lim_{\delta\downarrow0}
  \Big|&\int_0^T \psi(t) u^-_0(t)\,\d t-
  \iint_{\Omega} \psi(t)\varphi_\delta(s) \,\hat \mu_0(\d s\d t)\Big|\le
  \iint_{[0,T]\times\{-\kappa\}} \psi(t)\,\hat\sigma_0(\d s\d t)=0,
\end{align*}
since $\hat\sigma_0\ll\Lambda$.

On the other hand, recalling that $\hat\mu_0=\hat
m_0\,\Lambda$ and $\hat m_0\in
L^1_\lambda(0,T;W^{1,1}(\rmS))$,
an analogous argument yields for $\hat m_0^-(t):=
\hat m_0(t,-\kappa)$,
\begin{align*}
  \Big|\int_0^T \psi(t)\hat m^-_0(t)\,\dmm \lambda t-
  \iint_{\Omega} \psi(t)\varphi_\delta(s) \,\hat \mu_0(\d s\d t)\Big|
  &\le
  \iint_{\Omega} \psi(t)\varphi_\delta(s) \big|\hat m^-_0(t)-
  \hat m_0(t,s)\big|\,\d s\,\dmm \lambda t
  \\&\le
  \iint_{\Omega} \Phi_\delta(s) \big|\hat g_0(t,s)\big|\,\d s\,\dmm\lambda t
  \  \stackrel{\delta\downarrow0}{\longrightarrow}\ 0.
\end{align*}
Since $\psi$ is arbitrary, we conclude that 
\begin{equation}
\label{eq:mvsu}
\hat m_0^\pm\,\lambda= u_0^\pm\L^1.
\end{equation}

\smallskip
\textbf{Step 8: Passing to the limit in the continuity equation (statement \ref{th:continuityequation})}. This step is the same as in the proof of Theorem~\ref{th:compactness-bounded}.

\smallskip
\textbf{\boldmath Conclusion: Vanishing of the singular part of $\hat \nu_0$, i.e. $\hat w_0\lambda^\perp=0$.}
From~\eqref{eq:mvsu} it follows that $\hat
m_0^\pm(t)= 0$ for $\lambda^\perp$-a.e. $t\in [0,T]$.
On the other hand, \eqref{eq:3} yields for $\lambda^\perp$-a.e. $t\in
[0,T]$
and for every $\eta>0$ and $\hat s(t)$ with $\hat m_0(\hat s(t),t)>0$,
\begin{align*}
  +\infty&>\frac 12\int_\rmS \Big(\frac{\hat w_0(t)^2}{\hat m_0(t,s)}+
  \frac{\hat g_0(t,s)^2}{\hat m_0(t,s)}\Big)\,\d s\ge
  \frac 12\int_\rmS \Big(\frac{\hat w_0(t)^2}{\eta+\hat m_0(t,s)}+
  \frac{\partial_s \hat m_0(t,s)^2}{\eta+\hat m_0(t,s)}\Big)\,\d s
  \\&\ge
  |\hat w_0(t)|\int_\rmS |\partial_s \log(\eta+\hat m_0(t,s))|\,\d s\ge
  2|\hat w_0(t)|\Big|\log\eta-\log(\eta+\hat m_0(t,\hat s(t)))\Big|.
\end{align*}
Since $\eta>0$ is arbitrary, we conclude that $\hat w_0(t)=0$
$\lambda^\perp$-a.e.\ $t\in[0,T]$.
\end{proof}

The Lemma below is similar to many other duality results (see e.g.~\cite[\textsection2.6]{AmbrosioFuscoPallara00} or~\cite[Lemma~9.4.4]{AmbrosioGigliSavare05}) and seems to have some wider usefulness. We state it in $\R^d$ for generality.
\begin{lemma}
\label{lemma:supbdd}
Let $\Omega\subset\R^d$.
For  $\mu\in \M(\Omega)$ and $\nu\in\M(\Omega;\R^d)$,
\begin{equation}
\label{hyp:supbdd}
\frac12 \int_\Omega \Bigl|\frac{\d\nu}{\d\mu}\Bigr|^2 \, \d\mu = \sup \Biggl\{\int_\Omega \bigl[\,a\,\d \mu + b\cdot \d \nu\bigr] 
: \ a\in C_b(\Omega),\ b\in C_b(\Omega;\R^d),\  a+\frac{|b|^2}2\leq 0\ \Biggr\}.
\end{equation}
In particular, if the right-hand side is finite, then $\nu\ll\mu$ and
$\frac{\d\nu}{\d\mu} \in L^2_\mu(\Omega)$. 
\end{lemma}

\begin{proof}
  We write $\F(\nu|\mu)$ for the left-hand side, and $\F'(\nu|\mu)$
  for the right-hand side. We first show that $\F'(\nu|\mu)\leq
  \F(\nu|\mu)$. If $\nu$ is not absolutely continuous with respect to
  $\mu$, then $\F(\nu|\mu)=\infty$, and there is nothing to prove; if
  $\nu\ll\mu$, then we can write $\nu=f\mu$.  
%  and since the convex
%  function
%\begin{displaymath}
%[0,\infty)\times \R^d\ni (u,w) \mapsto
%\begin{cases}
%  \displaystyle\frac{|w|^2}{2u}&\text{if }u>0,\\
%  +\infty&\text{if }u=0,\ w\neq 0,\\
%  0&\text{if }u=0, \ w=0,
%\end{cases}
%\end{displaymath}
%is the convex conjugate of the indicator function of the set $\{(a,b)\in \R\times\R^d: a+ |b|^2/2\leq 0\}$, we have 
 For all $a$ and $b$ continuous, bounded, and satisfying $a+|b|^2/2\leq0$, we have
\[
\int_\Omega \bigl[a\,\d \mu + b\cdot \d \nu\bigr] 
= \int_\Omega \bigl[a + b\cdot f\bigr] \, d\mu
\leq \int_\Omega \Big[a + \frac{|b|^2}2 +\frac{|f|^2}{2}\Big] \, d\mu \leq \int_\Omega \frac{|f|^2}{2} \, d\mu = \F(\nu|\mu).
\]

To prove the opposite inequality, we assume that
$\F'(\nu|\mu)<\infty$, and first show that $\nu\ll\mu$. Suppose not;
then there exists a Borel set $A\subset\Omega$ such that $\mu(A)=0$
and $\nu(A)\not=0$. Take $c>0$, set $a = -c\chi_A$, and define a
sequence $a_n\in C_b(\Omega)$ such that $a_n\uparrow a$.  Then $\int
a_n\,\d\mu\to0$ as $n\to\infty$. On the other hand, setting $b_n :=
\sqrt{-2a_n} \nu(A)/|\nu(A)|$, we have $\int b_n\cdot \d \nu\to
\sqrt{2c}\, |\nu(A)|>0$. Since $c$ is arbitrary, this violates the
finiteness of $\F'(\nu|\mu)$, and therefore $\nu\ll\mu$.

Again writing $\nu = f\mu$, with $f\in L^2(\mu)^d$, we now choose
$b_n\in C_b(\Omega)$ such that $b_n\to f$ in $L^2(\mu)^d$, so that
$\int b_n\cdot f \,\d \mu\to \int |f|^2\, \d\mu =
2\F(\nu|\mu)$. Setting $a_n:= -|b_n|^2/2$ we have $a_n\to -|f|^2/2$ in
$L^1(\mu)$, and therefore $\int a_n\, \d\mu\to-\F(\nu|\mu)$. The
result follows.
\end{proof}

The above dual characterization \eqref{hyp:supbdd} of the property
$\frac{\d\nu}{\d\mu} \in L^2_\mu(\Omega)$ will now be used to
characterize the limits in Step 6 of the above proof. 

\begin{lemma}\label{prop:L2bdd}
If $u_n\longweakstarto \mu$ and $w_n\longweakstarto \nu$, and
\[
\sup_n \int_\Omega \frac{|w_n|^2}{u_n}\, \d x =:C <\infty,
\]
then $\nu \ll \mu$ with $\frac{\d\nu}{\d\mu} \in L^2_\mu(\Omega)$
and 
\begin{equation}
\label{lsc:w2u}
\int_\Omega \,\Bigl|\frac{\d\nu}{\d\mu}\Bigr|^2 \, \d\mu \leq \liminf_{n\to\infty}
\int_\Omega \frac{|w_n|^2}{u_n}.
\end{equation}
\end{lemma} 

\begin{proof} For each pair $(a,b) $ as in the right-hand side of 
\eqref{hyp:supbdd} we have 
\[
C\geq \int_\Omega \frac{|w_n|^2}{u_n}\, \d x \geq \int_\Omega \big[ a u_n + b\cdot w_n \big]\,\d x \to \int_\Omega \big[ 
a \,\d \mu + b\cdot \d \nu\big]. 
\]
Thus, the hypothesis of Lemma~\ref{lemma:supbdd} is satisfied and
$\frac{\d\nu}{\d\mu} \in L^2_\mu(\Omega)$ follows. 

Moreover, choosing a pair $(a,b)$ in \eqref{hyp:supbdd} that
approximates  the left-hand side in \eqref{lsc:w2u}, we also 
obtain the desired estimate \eqref{lsc:w2u}.
\end{proof}

\section{Lower bound}
\label{sec:lowerbound}

\begin{theorem}[Lower bound]
\label{th:inner_lower_bound}
{Under the same conditions as in}
Theorems~\ref{th:compactness-bounded} or~\ref{th:compactness}
 {let us assume,}
without loss of
generality,
 that $\rho_\e(t)\weakstarto \rho_0(t)$ for every $t\in [0,T]$.
Then
\begin{equation}
\label{ineq:lower-bound}
\J_0(\rho_0)\leq \liminf_{\e\to0}
\J_\e(\rho_\e)\quad\text{and}\quad
\E_0(\rho_0(t))\le
\liminf_{\e\to0}\E_\e(\rho_\e(t))\quad
\text{for every }t\in [0,T].
\end{equation}
\end{theorem}

\begin{proof}
The lower semicontinuity of the entropy functionals under weak
convergence
is well known, see e.g.~\cite[Lemma~9.4.3]{AmbrosioGigliSavare05} or~\cite[Lemma~6.2]{AmbrosioSavareZambotti09}.

Turning to $\J_\e$, 
we can suppose by Theorem~\ref{th:compactness} (which contains Theorem~\ref{th:compactness-bounded})
that 
\[\hat u_\e \longweakstarto \hat \mu_0 = \hat m_0\Lambda,
\qquad
\partial_s \hat u_\e \longweakstarto \hat\varsigma_0 = \hat g_0\Lambda,
\qquad\text{and}\qquad
\hat w_\e \longweakstarto \hat \nu_0 = \hat n_0 \Lambda.
\]
Setting $\hat u_0 = \hat m_0\ell$ as in the proof of Theorem~\ref{th:compactness}, we also have $\hat g_0 \ell= \partial_s \hat u_0$. By~\eqref{lsc:w2u} we then have
\[
\frac12\int_0^T\!\!\int_\rmS \;\biggl[\Bigl(\frac{\d\hat\nu_0}{\d\hat\mu_0}\Bigr)^2
+ \Bigl(\frac{\d\hat\varsigma_0}{\d\hat\mu_0}\Bigr)^2\biggr]\, \d\hat \mu_0
\leq \liminf_{\e\to0} \J_\e(\rho_\e).
\]
We now discard the singular part $\hat m_0\Lambda^\perp$ of $\hat\mu_0$ and again write $\hat w_0 := \hat n_0\ell$, by which we find
\[
\frac12\int_0^T\!\! \int_\rmS \Big(\frac{\hat w_0^2}{\huz}+
  \frac{ \partial_s \hat u_0^2}{\huz}\Big)\, \d s\,\d t\le
\frac12\int_0^T\!\!\int_\rmS \;\biggl[\Bigl(\frac{\d\hat\nu_0}{\d\hat\mu_0}\Bigr)^2
+ \Bigl(\frac{\d\hat\varsigma_0}{\d\hat\mu_0}\Bigr)^2\biggr]\, \hat u_0\,\d s\d t
\leq \liminf_{\e\to0} \J_\e(\rho_\e).
\]
Recalling that the traces of $\hat u_0$ at $s=\pm \kappa$ coincide
with $u^\pm_0$, and that $\hat w_0(t,s)=\hat w_0(t)$ is constant with
respect to $s$ with $\hat w_0=\frac 12\dot u^+$, we see that
for a.e.~$t$ the integrand in the left-hand side of the previous
inequality satisfies
\begin{displaymath}
\Md{\hat w_0(t)}{u_0^\pm(t)} \leq 
  \frac12 \int_\rmS \Big(\frac{\hat w_0(t)^2}{\huz(t,s)}+
  \frac{ \partial_s \hat u_0(t,s)^2}{\huz(t,s)}\Big)\, \d s  .
\end{displaymath}
This implies the lower bound on $\J_\e$ and concludes the proof of Theorem~\ref{th:inner_lower_bound}.
\end{proof}

\section{The minimization problem defining $M$ and interpolation}
\label{sec:M}

The minimization problem defining $M$ is
\begin{equation}
\label{def:M_secmin}
\Md w{u^\pm} := \inf_u \left\{ \;\frac12\int_\rmS \Bigl[ \frac{w^2}{u(s)} + \frac{
 u'(s)^2 } {u(s)}\Bigr]\, \d s\;:\quad
u(\pm\kappa) = u^\pm\;\right\}.
\end{equation}
This minimization problem gives rise to a natural \emph{interpolation} of the boundary values $u^\pm$, which we study in the following theorem.

\begin{theorem}
\label{th:propsM}
Let $u[w,u^\pm](\cdot)$ be the solution of the minimization problem $M(w,u^\pm)$. Then the mapping
\[
(w,u^\pm)\mapsto u[w,u^\pm]
\]
is well-defined and continuous from $\R\times (0,\infty)^2$ into $C^2(\rmS)$.
The function $(w,u^\pm) \mapsto M(w,u^\pm)$ is convex,  smooth away from $u^\pm = 0$, minimal at $w=0$ and $u^+=u^-$, and satisfies $M(w,u^\pm) = M(w,u^\mp)$.

If $u^\pm \in C^2([0,T]; [\delta,\infty))$ for some $\delta>0$, then the function
\begin{equation}
\label{mapping:a-time-dep}
(t,s) \mapsto u\Bigl[\frac12 \dot u^+(t),u^\pm(t)\Bigr](s) 
\end{equation}
is an element of $C^1([0,T]\times \rmS)$.
\end{theorem}

\begin{proof}
By the transformation $z = \sqrt u$ we can rewrite the minimization problem~\eqref{def:M} as
\[
\inf_z \int_\rmS \Bigl[\frac{w^2}{2z(s)^2} + 2z'(s)^2\Bigr]\, \d s
\; : \quad z(\pm\kappa) = \sqrt {u^\pm}.
\]
The corresponding stationarity equation is 
\begin{equation}
\label{ELeq:w}
-z'' - \frac{w^2}{4z^3} = 0, \qquad \quad z(\pm\kappa) = \sqrt {u^\pm},
\end{equation}
which implies that any solution $z$ is concave and therefore $z\geq \min \sqrt{u^\pm}$,
or $u\geq \min u^\pm$.

Since $u^\pm>0$, the existence and uniqueness of the solution $u$ of~\eqref{def:M_secmin}, or equivalently of the solution $z$ of~\eqref{ELeq:w}, are classical, and the continuity follows from classical results for the continuous dependence of the solutions of elliptic problems on parameters. Similarly, if $u^\pm$ is a function $u^\pm\in C^2([0,T])^2$ and bounded away from zero, then the solution $u[\dot u^+(\cdot)/2,u^\pm(\cdot)]$ is $C^1$ on its domain (note that one degree of differentiation in time is lost since $\dot u^+ $ appears as a parameter in the equation for $z$).

\medskip
The symmetry and minimality properties of $M$ are immediate. To prove the convexity of $M$, take $(w_1,u_1^\pm)$ and $(w_2,u^\pm_2)$ with $M(w_1,u_1^\pm)$, $M(w_2,u_2^\pm)<\infty$, $\lambda\in[0,1]$, and let $u_1$ and $u_2$ be the corresponding minimizers. Since $(u,w)\mapsto w^2/u$ is convex, it follows that 
\[
\frac{(\lambda w_2 + (1-\lambda)w_1)^2}{8(\lambda u_2(s) + (1-\lambda) u_1(s))}
\leq \lambda \frac{w_2^2}{8u_2(s)} + (1-\lambda) \frac{w_1^2}{8u_1(s)},
\]
with a similar inequality for the second term in~\eqref{def:M_secmin}. Since $\lambda u_2 + (1-\lambda)u_1$ is admissible for $M(\lambda w_2 + (1-\lambda)w_1, \lambda u^\pm_2 + (1-\lambda)u^\pm_1)$, we then have
\begin{align*}
M(\lambda w_2 &+ (1-\lambda)w_1, \lambda u^\pm_2 + (1-\lambda)u^\pm_1)\leq \\
&\leq \int_0^1 \Bigl[\frac{(\lambda w_2 + (1-\lambda)w_1)^2}{8(\lambda u_2(s) + (1-\lambda) u_1(s)}
+ \frac{(\lambda u_2' + (1-\lambda)u_1')^2}{2(\lambda u_2(s) + (1-\lambda) u_1(s)}\Bigr]\, ds\\
&\leq \lambda \int_0^1 \Bigl[\frac{w_2^2}{8u_2(s)} + \frac{u_2'(s)^2}{2u_2(s)}\Bigr]\, ds
+ (1-\lambda) \int_0^1 \Bigl[\frac{w_1^2}{8u_1(s)} + \frac{u_1'(s)^2}{2u_1(s)}\Bigr]\, ds\\
&= \lambda M(w_2,u_2^\pm) + (1-\lambda) M(w_1,u_1^\pm).
\end{align*}
This concludes the proof of Theorem~\ref{th:propsM}.
\end{proof}

As indicated at the end of the introduction, we can find good lower and upper bounds on the integrand~$M$, which are given in the following theorem.

\begin{theorem}\label{th:M.est} For all $u^\pm>0$ and all
 $w\in \R$ we have the estimate 
 \begin{equation}
 \label{eq:M.est}
   w(\log u^+-\log u^-) 
%   + \frac{\bigl(4\kappa^2w^2-(u^+-u^-)\bigr)^2}{16\kappa^2(u^++u^-)}
   \leq M(w;u^\pm) \leq
   \frac{\log u^+-\log u^-}{4\kappa(u^+-u^-)} \big( 4\kappa^2w^2 + (u^+-u^-)^2\big),
 \end{equation}
where both inequalities are equalities if and only if $w = (u^+-u^-)/2\kappa$. 
In this case the minimizer~$u$ in the definition
\eqref{def:M_secmin}  of
$M(w;u^\pm)$ is the affine interpolation $u(s)= (\kappa+s)u^+/2\kappa + (\kappa-s)u^-/2\kappa$.
\end{theorem}

\begin{remark}
Note that the left-hand side of~\eqref{eq:M.est} can be interpreted as $2\langle\E'(u^\pm),w\rangle$ (see Section~\ref{subsec:structureJ0}), implying that $M(w;u^\pm)\geq 2\langle\E'(u^\pm),w\rangle$ and therefore $\J_0(u^\pm)\geq0$ for all $u^\pm$. 
%In addition, if $\hat J_0(\huz)=0$, then for almost all $t$ the function $\huz(t,\cdot)$ achieves the minimum of $M(\huz(t),\dhuz(t))$.
\end{remark}

\begin{proof} We define the functional $J(w;u)=\frac12\int_\rmS
 \frac1u(w^2 + {u'}^2)\,ds$ such that $M$ is obtained by
 minimizing $J(w;u)$ over all $u$ satisfying the boundary conditions $u(\pm\kappa)=u^\pm$. 

The lower estimate follows by neglecting the nonnegative term in
\[
J(w;u)\;=\;\int_\rmS \frac1{2u}\big(w-u'\big)^2 + w \int_\rmS
\frac {u'}u \;\geq\; w(\log u(\kappa)-\log u(-\kappa))
\]
and using the boundary conditions. We also see that equality holds if
and only if $u' \equiv v/2$, which implies $v= 4(m{-}a)$. 

The upper estimate is obtained by testing with the affine function
$u(s) = (\kappa+s)u^+/2\kappa + (\kappa-s)u^-/2\kappa$. Obviously, the lower estimate and the upper
estimate coincide for $w=(u^+-u^-)/2\kappa$. Hence, the result is
proved. 
\end{proof} 

%The fact that optimality occurs at affine functions is reminiscent of
%the fact that the limits $\hat u_0(t,s)$ are affine in $s\in [0,1]$,
%see Theorem \ref{th:compactness_solutions}.

The fact that optimality occurs at affine functions also gives a characterization of the limit $\huz$ of a sequence of \emph{solutions} $\hue$:
\begin{theorem}
\label{th:compactness_solutions}
Let $\rho_\e$ be a sequence of solutions of~\eqref{pb:main} such that $\E_\e(\rho_\e(0))$ converges as $\e\to0$. Then the assertions of Theorem~\ref{th:compactness} hold, and in addition $\huz$ is affine in $s$:
\[
\text{for almost all }t,s, \qquad \huz(t,s) = \frac{\kappa+s}{2\kappa} \huz(t,\kappa) + \frac{\kappa-s}{2\kappa} \huz(t,-\kappa).
\]
\end{theorem}

\begin{proof}
%By the maximum principle for parabolic equations~\cite{ProtterWeinberger67} the uniform bound on the initial data $u_\e(0,\cdot)$ implies a uniform $L^\infty$-bound on the solution $u_\e$ for all time. Since solutions of~\eqref{pb:main_u} are also smooth, the conditions of Theorem~\ref{th:compactness-bounded} are satisfied. 
%
The transformed solutions $\hue$ satisfy the equation
\[
\hgdens \partial_t \hue = \partial_{ss}\hue.
\]
The density $\hgdens$ concentrates on to the boundary points $s=\pm\kappa$, implying that in the interior of the interval $\rmS$ the equation formally reduces to $0=\partial_{ss} \hue$. Using classical methods for partial differential equations one can convert this observation into a proof that the limit $\huz$ is affine for each~$t$. 

Instead we prefer to stay within the realm of the gradient-flow structure. 
Since the $\rho_\e$ are solutions, $\A_\e(\rho_\e)=0$; by Theorem~\ref{th:inner_lower_bound} and the assumption of convergence of the initial energies $\E_\e(\rho_\e(0))$, we have $\A_0(\rho_0)\leq 0$. Since $\A_0$ satisfies condition~\eqref{cond:nonnegativeA}, $\A_0(\rho_0)=0$. This  implies  that $\huz$ is a minimizer of $M$ for almost all $t$, and by Theorem~\ref{th:M.est} it is therefore affine for almost all~$t$.
\end{proof}

\section{Recovery sequence}
\label{sec:upperbound}

\begin{theorem}[Recovery sequence]
\label{th:recovery_sequence}
Let $u^\pm \in AC(0,T; \R)$ be such that $\J_0(u^\pm)<\infty$. Then there exists a sequence $\hue\in C^1([0,T]\times\rmS)$ such that $\hue(\cdot,\pm\kappa) \to u^\pm$ in $L^1(0,T)$, $\hE_\e(\hue(0))\to\E_0(u^\pm(0))$, $\hE_\e(\hue(T))\to\E_0(u^\pm(T))$, and $\hJ_\e(\hue)\to \J_0(u^\pm)$.
\end{theorem}

\begin{remark}
By this result the sequence $\hue$ and its other forms $u_\e$, $\rho_\e$, and $\hat \rho_\e$ converge in the different senses provided by Theorem~\ref{th:compactness-bounded}.
\end{remark}

%\begin{remark}
%Since the sequence constructed in Theorem~\ref{th:recovery_sequence} satisfies the conditions of Theorem~\ref{th:compactness}, we can also choose it such that it converges in the senses described by Theorem~\ref{th:compactness}.
%\end{remark}

\begin{proof}
By a diagonal argument, and using the lower bound~\eqref{ineq:lower-bound}, it is sufficient to prove the following approximation result: given $\delta>0$, there exists a sequence $(\hue^\delta)_{\e>0}$ with $\hue^\delta\hge\in C^1([0,T]\times\rmS)$ and
\begin{equation}
\label{recovery:approximation}
\limsup_{\e\to0} \  \max\Bigl\{\  
\|\hue^\delta(\cdot,\pm\kappa)-u^\pm\|_{L^1(0,T)^2}, \ 
\bigl|\hE_\e(\hue^\delta(T))-\E_0(u^\pm(T))\bigr|,\ \\ 
\bigl|\hE_\e(\hue^\delta(0))-\E_0(u^\pm(0))\bigr|\ \Bigr\} \leq \delta,
\end{equation}
and
\begin{equation}
\label{recovery:approximation'}
\limsup_{\e\to0}   
\hJ_\e(\hue^\delta)-\J_0(u^\pm)\leq \delta.
\end{equation}

%\begin{equation}
%\label{recovery:approximation}
%\limsup_{\e\to0}\|\hue^\delta(\cdot,\pm\kappa)-u^\pm\|_{L^1(0,T)} \leq \delta \qquad\text{and}\qquad
%\limsup_{\e\to0}|\hJ_\e(\hue^\delta)-\J_0(u^\pm)| \leq \delta.
%\end{equation}
We now prove this approximation result in several steps.
First note that $u^\pm\in W^{1,1}(0,T)$, and that the finiteness of $\J_0(u^\pm)$ implies that $u^++u^-$ is constant in time (say $2\mass$) and therefore $0\leq u^\pm\leq 2\mass$ and $\dot u^+=-\dot u^-$. 
To simplify we only specify the value $u^-$ at $-\kappa$, and consider the corresponding value $u^+$ at $+\kappa$ as defined by the condition of constant mass.

We first approximate $u^-$ by a function that is bounded away from
zero and from $2\mass$. We do this by setting $y_\eta^- := \mass +
(1-\eta)(u^--\mass)$, for some small $\eta$; as $\eta\to0$,
$y_\eta^-\to u^-$ in $W^{1,1}(0,T)$. The function $y_\eta^-$ is bounded away
from $0$ and $2\mass$; the convexity of $M$ and the fact that it vanishes when $w=0$ and $u^+=u^-$ (Theorem~\ref{th:propsM}) imply that for almost all $t$,
$M(\dot y_\eta^+(t)/2;y^\pm_\eta(t))$ is decreasing in~$\eta$, and that
$M(\dot y^+_\eta(t)/2; y^\pm_\eta(t))\uparrow M(\dot u^+(t)/2;u^\pm(t))$ as
$\eta\downarrow0$. This implies that
\[
\int_0^T M\Bigl(\frac12 \dot y_\eta(t)^+; y^\pm_\eta(t)\Bigr)\, \d t \longrightarrow \int_0^T M\Bigl(\frac12\dot u^+(t);u^\pm (t)\Bigr)\, \d t 
\qquad\text{as}\ \eta\to 0.
\]
Similarly $\E_0(y_\eta^\pm(T))\to\E_0(u^\pm(T))$ and $\E_0(y_\eta^\pm(0))\to\E_0(u^\pm(0))$,  implying that for given $\delta>0$ we may choose $\eta>0$ such that 
\begin{multline}
\label{th:recovery:est1}
\max\Bigl\{\ \|y_\eta^\pm-u^\pm\|_{L^1(0,T)^2},\ \bigl|\E_0(y_\eta^\pm(T))-\E_0(u^\pm(T))\bigr|,
 \\
\bigl|\E_0(y_\eta^\pm(0))-\E_0(u^\pm(0))\bigr|,\ 
\bigl|\J_0(y_\eta^\pm)-\J_0(u^\pm)\bigr|\ \Bigr\}\leq \frac\delta4.
\end{multline}
We fix this number $\eta$.

The next step is to smoothen $y_\eta$. We approximate $y_\eta$ in $W^{1,1}(0,T)$ by convolution to give a $\widetilde y\in C^2([0,T])$, while preserving the pointwise upper and lower bounds. Because $M$ is convex, it follows that 
\begin{equation}
\label{th:recovery:est2}
\J_0(\widetilde y^\pm)\leq \J_0(y_\eta^\pm),
\end{equation}
and we can choose $\widetilde y$ such that 
\begin{equation}
\label{th:recovery:est2'}
\max\Bigl\{\ \|\widetilde y^\pm -y_\eta^\pm\|_{L^1(0,T)^2},\ 
\bigl|\E_0(\widetilde y^\pm(T))-\E_0(y_\eta^\pm(T))\bigr|,\ 
\bigl|\E_0(\widetilde y^\pm(0))-\E_0(y_\eta^\pm(0))\bigr|,\ \Bigr\}\leq \frac\delta4.
\end{equation}
%\begin{equation}
%\label{th:recovery:est2}
%\|\widetilde y-\huz\|_{L^2(0,T;L^2(0,1))} \leq \frac\delta2 \qquad\text{and}\qquad
%|\hJ_0(\widetilde y)-\J_0(\huz)| \leq \frac\delta4.
%\end{equation}

We now interpolate $\widetilde y$ by Theorem~\ref{th:propsM} without changing notation; note that then $\widetilde y\in C^1([0,T]\times \rmS)$.
Since $\widetilde y$ is fixed and $C^1$, it follows that as $\e\to0$, the corresponding function $\hat w_\e$, defined by $\partial_t (\widetilde y\hge) + \partial s \hat w_\e=0$,  is uniformly bounded and satisfies
\[
\forall s\in[-\kappa,\kappa), \qquad \hat w_\e(t,s)  = \int_{-\kappa}^s \dot{\widetilde y}(t,\sigma) \, \hge(\d\sigma)
\to \frac12 \dot{\widetilde y}(t,-\kappa) = \frac12\dot {\widetilde y}^-(t). 
\]
Therefore
\[
\lim_{\e\to0} \int_0^T\!\!\int_\rmS \frac{\hat w_\e^2}{\widetilde y}\, \d s\d t = \int_0^T\!\!\int_\rmS \frac{w^2}{\widetilde y}\, \d s\d t \qquad
\text{with }w(t) = \frac12 \dot {\widetilde y}^-(t),
\]
and we have for sufficiently small $\e>0$ that
\begin{equation}
\label{th:recovery:est3}
|\hJ_\e(\widetilde y) - \J_0(\widetilde y^\pm)| \leq \frac\delta 4.
\end{equation}
Similarly, since $\hge\weakstarto \hgz$, we have for sufficiently small $\e>0$ that 
\begin{equation}
\label{th:recovery:est3b}
\max\Bigl\{ \ \bigl|\hE_\e(\widetilde y(T)) - \E_0(\widetilde y^\pm(T))\bigr|, \ 
  \bigl|\hE_\e(\widetilde y(0)) - \E_0(\widetilde y^\pm(0))\bigr| \ \Bigr\} \leq \frac\delta4.
\end{equation}

The final step is to approximate $\tilde y$ by a function of the right mass. Since $\widetilde y\in C^2([0,T]\times\rmS)$, the mass discrepancy
\[
\widetilde m(t) := \int_\rmS \widetilde y(t,s)\, \dmm \hge s - \int_\rmS \widetilde y (t,s) \,\dmm \hgz  s
\]
converges to zero uniformly on $[0,T]$. Setting 
\[
\hue^\delta(t,s) := \widetilde y(t,s) - \widetilde m(t),
\]
we find that for sufficiently small $\e$
%Choose $\varphi\in C_c^\infty((-\kappa,\kappa))$ with $\varphi\geq0$, $\varphi\not\equiv0$, and choose the sequence $a_\e\in\R$ such that 
%\[
%\hue^\delta(t,s) := \widetilde y(t,s) + a_\e(t)\varphi(s)
%\]
%satisfies $\int \hue^\delta\, d\hat \gamma_\e = m$. Since $\int \widetilde y\, d\hat \gamma_\e\stackrel{\e\to0}\longrightarrow \int \widetilde y \,d\hat \gamma_0 = m$, this implies $a_\e\varphi\hge \to0$ in $L^1(\rmS)$ for all $t$. For sufficiently small~$\e$, we therefore have
\begin{multline}
\label{th:recovery:est4}
\max\Bigl\{\ \|\hue^\delta(\cdot,\pm\kappa) - \widetilde y^\pm\|_{L^1(0,T)^2},\ 
\bigl|\hE_\e(\hue^\delta(T))-\hE_\e(\widetilde y(T))\bigr|,
 \\
\bigl|\hE_\e(\hue^\delta(0))-\hE_\e(\widetilde y(0))\bigr|,\ 
\bigl|\hJ_\e(\hue^\delta)-\hJ_\e(\widetilde y)\bigr|\ \Bigr\}\leq \frac\delta4.
\end{multline}
%\begin{equation}
%\label{th:recovery:est4}
%\|\hue^\delta - \widetilde y \|_{L^2(0,T;L^2(0,1))} \leq \frac\delta2 
%\qquad\text{and}\qquad
%|\hJ_\e(\hue^\delta) - \hJ_\e(\widetilde y)| \leq \frac\delta4,
%\end{equation}
The claims~\eqref{recovery:approximation} and~\eqref{recovery:approximation'} then follow from combining the estimates~\eqref{th:recovery:est1}, \eqref{th:recovery:est2}, \eqref{th:recovery:est2'}, \eqref{th:recovery:est3}, \eqref{th:recovery:est3b}, and~\eqref{th:recovery:est4}.
\end{proof}

\section{Connections with stochastic particle systems}
\label{sec:ldp}

The mathematical results of this paper make important use of Definition~\ref{def:GF} of a gradient-flow solution. This formulation is more than a mathematical convenience:  it arises naturally when considering equation~\eqref{pb:main} as the deterministic limit of a stochastic system of particles. We now explain this connection and its consequences. 

\medskip
Fix $\e$ for the moment. Consider a collection of $n$ independent particles, each of which performs a Brownian motion in a potential landscape given by the energy function $H/\e$. 
Equation~\eqref{pb:main} is the \emph{continuum} or \emph{hydrodynamic limit} of this system of particles, as the number $n$ of particles tends to infinity. One way of describing this limit is by considering the \emph{empirical measure}
\[
L_n: [0,T]\to\M([-1,1]), \qquad
L_n(t) := \frac1n \sum_{i=1}^n \delta_{\xi_i(t)},
\]
where $\xi_i(t)$ is the position in $[-1,1]$ of particle $i$ at time $t$. As $n\to\infty$, with probability $1$ this empirical measure converges weakly to a limit measure $\rho(t)$ at every time $t$, and this limit measure solves  the equation~\eqref{pb:main}.\footnote{We deliberately disregard the role of initial conditions for the moment; one could, for instance, choose a single delta function $\delta_a$ as the initial datum for~\eqref{pb:main}, and have all particles start at $\xi=a$.}

Given this connection, a \emph{large-deviation} result characterizes the probability of finding the empirical measure $L_n(t)$ \emph{far} from the solution $\rho$ of~\eqref{pb:main}. Such a result roughly takes the form
\begin{equation}
\label{ldp}
P_n (L_n \approx \tilde \rho) \sim \exp({-}n I(\tilde \rho)),
\end{equation}
in terms of a \emph{rate functional} $I$. A rigorous version of this statement can be found, for instance, in~\cite[Th. 13.37]{FengKurtz06}.

The surprising feature, however, is that for this system of particles, the rate functional $I$ above is \emph{exactly} equal to the functional $\A_\e$ in~\eqref{def:AWass} (see e.g.~\cite{KipnisOlla90} or~\cite[Th. 13.37]{FengKurtz06}). This feature has several consequences.
\begin{enumerate}
\item Since $\A_\e=I$, the large-deviation result~\eqref{ldp} gives an
 alternative explanation why $\A_\e\geq0$ and why $\A_\e=0$ implies a
 solution of the deterministic system. The positivity of $I$, and
 therefore of $\A_\e$, arises directly from the property~\eqref{ldp} and
 the fact that probabilities are bounded by $1$; and since the
 hydrodynamic limit is assumed with probability $1$, the solution~$\rho$ of the limit equation~\eqref{pb:main} necessarily satisfies
 $\A_\e(\rho)=I(\rho) =0$.\footnote{A careful treatment of this argument
   actually requires a more precise definition of~\eqref{ldp} and a
   discussion of topology; we omit both.}
\item In Section~\ref{subsec:differentstructures} we mentioned that
 there exist at least two different gradient-flow structures for
 equation~\eqref{pb:main}. The fact that one of these structures
 arises in the large-deviation description of this stochastic system,
 may be interpreted to signify that this gradient-flow structure is
 more `natural'---at least when we view~\eqref{pb:main} as arising
 from this specific stochastic particle system. Of course, there may
 well be a different stochastic system whose large-deviation
 behaviour is related to the structure~\eqref{eq:GS_lin}, and there
 may be other arguments that favour other structures.
\item This connection provides an answer to the question, often heard,
 `why does the Wasserstein metric figure in this gradient-flow
 structure?', since the Wasserstein dissipation arises directly from
 the large-deviation behaviour. However, a complete answer requires
 describing the large-deviation result in some detail, which would
 take us too far; see~\cite{AdamsDirrPeletierZimmer10TR} for a
 detailed discussion.
\item In the context of a large-deviation result, it is natural to
 consider sets of the form $\{\tilde \rho: I(\tilde\rho)\leq
 \delta\}$ for $\delta>0$; these correspond to collections of `least
 unlikely' states, in the sense that their probability vanishes no
 faster than $\rme^{-n\delta}$. Again, this ties in with the results
 proved above, in which we do not assume $\A_\e=0$, but only
 boundedness of $\E_\e(\rho_\e(0))$ and~$\J_\e(\rho_\e)$.
\end{enumerate}

The connection with large-deviation principles also explains the structure of $M$ in~\eqref{def:M}. The well-known \emph{contraction principle} describes how rate functions transform under projection, i.e. under loss of information. Suppose that $I$ is a rate function describing the behaviour of a sequence of probability measures $P_n$ on a space $X$, in the sense of~\eqref{ldp}. Let $p:X\to Y$ be a continuous map, and $Q_n := P_n\circ p^{-1}$ the corresponding probability measures describing the behaviour of the system after projection under $p$ onto $Y$. Then $Q_n$ satisfies the large-deviation principle~\cite[Th. III.20]{DenHollander00}
\[
Q_n(y) \sim \exp({-}n I_p(y)) \qquad
\text{with}\qquad
I_p(y) := \inf_{x\in X: p(x) = y} I(x).
\]
The form of the function $M$ can be understood in terms of this contraction principle. In the limit $\e=0$, the only information about $u_\e$ or $\hue$ that survives are the boundary values. Consequently the large-deviation behaviour of the system in the limit follows from the contraction principle by interpolating between the boundary values, in such a way as to minimize the functional over all missing information. The function $M$ is the direct consequence of this.

\section{Discussion}
\label{sec:discussion}

\emph{Passing to the limit in gradient flows.}
The aim of this work is to explore 
the potential of the Wasserstein gradient-flow structure of~\eqref{pb:main} for rigorous passing to the limit. For this specific system, we have succeeded to a large degree, and we comment below on the specific assumptions that we have made.

The property $\A=0$ is a reformulation of the concept of a \emph{curve of maximal slope}, which was introduced by DeGiorgi and co-authors (see e.g.~\cite{DeGiorgiMarinoTosques80}) as a metric-space generalization of a gradient flow. Sandier, Serfaty, and Stefanelli~\cite{SandierSerfaty04,Stefanelli08} appear to be the first to explore in detail the use of this structure for passing to the limit. 
Serfaty~\cite{Serfaty09TR} discusses the case of metric spaces, with obvious applications for the case of the Wasserstein metric. She leaves aside the question of compactness, however, and one of the main contributions of this paper is to show that appropriate compactness `in time' can also be obtained from the Wasserstein structure.

A related result is that of Ambrosio, Savar\'e, and Zambotti~\cite{AmbrosioSavareZambotti09}, who study entropy-Wasserstein gradient flows in a Hilbert space with respect to a weakly converging sequence of reference measures. Their approach first proves convergence of time-discrete approximations for fixed time, and then uses error bounds to prove convergence of the time-continuous solutions. 

\medskip

\emph{Assumptions.} 
The assumptions in the main theorems are the boundedness of the initial energy and of the dissipation function $\J_\e$, both of which are natural objects in the Wasserstein gradient flow. The relaxation of the condition $\A_\e=0$ to
the condition
$\sup_\e 
\J_\e< \infty$ is a broadening of scope: it implies that the
compactness result holds not only for solutions, with their
accompanying higher-regularity properties, but for a much wider class
of sequences. In addition, this class arises naturally in the context
of large deviations for an underlying stochastic particle system (see Section~\ref{sec:ldp}). 

However, a central tool is the mapping $\xi\mapsto s$, which desingularizes the diffusion term and allows for a more detailed study of the limit behaviour. This mapping is very specific for this problem, and it is an interesting question how to generalize it to singular systems described by different PDEs (e.g.\ higher-order~\cite{Otto98a,GiacomelliOtto01,Glasner03} or nonlocal parabolic equations~\cite{CarrilloMcCannVillani03,CarrilloMcCannVillani06}) or more complicated geometric spatial structure.

\medskip

\emph{Weak formulations and compactness.}
A related question arose during the work presented here: can our definition of a gradient-flow solution, Definition~\ref{def:GF}, be viewed as a \emph{weak form} of the gradient-flow equation~\eqref{eq:GradFlow}, 
\begin{equation}
\tag{\ref{eq:GradFlow}}
\dot z = -\nabla_G \E(z)?
\end{equation}
The straightforward answer to this question seems to be negative, since traditionally weak formulations serve to reduce regularity requirements, and in both cases the function $z$ is necessarily differentiable. Therefore shifting from~\eqref{eq:GradFlow} to Definition~\ref{def:GF} brings no advantage on that front. 

However, we argue here that a different aspect is just as important: the compactness and convergence properties of the formulation. As we have shown, solutions of a sequence of problems are compact in an appropriate way, and a subsequence converges to a limiting object that can be considered an appropriate generalization  of a gradient flow (see Definition~\ref{def:GF}, and also the discussion in Section~\ref{sec:ldp}). 

This ties in with the strongly related work of Herrmann and Niethammer~\cite{HerrmannNiethammer11}, that we discuss separately below. One aspect that this paper and~\cite{HerrmannNiethammer11} have in common is the reformulation of a nonlinear, singular differential equation as a parameter-dependent variational problem, thereby opening the door to methods of variational calculus.

\medskip

\emph{Choice of convergence.} 
We prove Gamma-convergence of the sequence of functionals $\J_\e$. 
If one is only interested in convergence of solutions, then this is actually too strong: it suffices to prove the lower bound inequality, Theorem~\ref{th:inner_lower_bound}. We prove the recovery sequence, Theorem~\ref{th:recovery_sequence} nonetheless, especially since it completes the picture of the convergence on $\J_\e$. 

Incidentally, the fact that Gamma-convergence is a natural form of convergence for large-deviation rate functionals (see Section~\ref{sec:ldp}) has also been recognized in the probabilistic literature~\cite{Leonard07TR}.

\medskip
\emph{The micro-problem.}
The transformation to the new spatial variable $s$ has the effect of blowing up the region in which the derivative of $u_\e$ is large. The resulting function $\hue$ is more regular, as reflected by the $H^1$-bound on $\surd \hue$ (see Figure~\ref{fig:transformation} and Remark~\ref{rem:desingularizing}). This blow-up argument is reminiscent of the `cell problem' in homogenization~\cite{Hornung97} or the `inner' and `outer layers' in singular perturbation theory~\cite{Verhulst05}. Similarly, the parallel convergence results of those two theorems reflect these separate behaviours at two different scales.

It is interesting to note, however, that for the lower-bound inequality one does not need to know much about the function $M$; actually, only its definition. Some additional information (the inequalities of Theorem~\ref{th:M.est}) is necessary to identify solutions of $\A_0=0$ as solutions of a corresponding differential equation. Other additional information (Theorem~\ref{th:propsM}) is necessary for the recovery sequence, Theorem~\ref{th:recovery_sequence}.

\medskip

\section{Comparison with a paper by Herrmann and Niethammer~\cite{HerrmannNiethammer11}}
\label{sec:HerrmannNiethammer}

This is not the first paper to give an answer to the question that was raised in~\cite{PeletierSavareVeneroni10}, \emph{Can we prove convergence using the Wasserstein gradient flow?} In~\cite{HerrmannNiethammer11}, Herrmann and Niethammer give a different (but again affirmative) answer. Here we briefly describe their approach and comment on the differences.

The authors of~\cite{HerrmannNiethammer11} build upon a solution concept for gradient flows based on an integrated form of the Rayleigh principle. This concept has been used before in nearly-finite-dimensional situations~\cite{NiethammerOtto01,NiethammerOshita10}, but its application in a truly parabolic context appears to be new. We describe it here in the case of a linear space $\Z$; the generalization to a manifold is straightforward. Given an energy functional $\E$ on $\Z$ and a Riemannian metric $G$ (see the introduction), it is  a straightforward observation that if $z$ is a solution of the gradient-flow equation~\eqref{eq:GradFlow}, then its time derivative $\dot z(t)$ at time $t$ is a minimizer of 
\begin{equation}
\label{def:RP}
v\mapsto \langle v,G(z(t))v\rangle + \langle\mathrm D\E(z(t)),v\rangle .
\end{equation}
Inspired by this, the authors of~\cite{HerrmannNiethammer11} define 
an \emph{integrated Rayleigh principle} as follows: an absolutely continuous function $z:[0,T]\to \Z$ satisfies this principle if its time derivative $\dot z$ minimizes the functional
\begin{equation}
\label{def:intRP}
\int_0^T \Bigl(\frac12 \langle v(t),G(z(t))v(t)\rangle
+ \bigl\langle\mathrm{D}\E(z(t)), v(t)\bigr\rangle\Bigr)\, \d t,
\end{equation}
among all $v:[0,T]\to\Z$.

In~\cite{HerrmannNiethammer11}, the authors first remark that at finite $\e>0$, the solution $u_\e$ of~\eqref{pb:main_u} is a minimizer of this integrated Rayleigh principle. In addition, the \emph{a priori} estimates of~\cite{PeletierSavareVeneroni10} provide appropriate compactness of the sequence $u_\e$. The central result is then that the integrated Rayleigh principle for the limiting function $u_0$ can be derived from the same principle for the solutions $u_\e$. 

\medskip

The work by Herrmann and Niethammer is interesting for various reasons. First, this solution concept merits to be considered more closely, and we make some comments on this below. Next, the authors themselves state as a drawback that the compactness that they use does not derive from the Wasserstein gradient-flow structure, but from the linear semi-group structure used in~\cite{PeletierSavareVeneroni10}. They pose in turn the question whether the compactness can be derived from the Wasserstein structure. And finally, what is exactly the relationship between the solution concepts of~\cite{HerrmannNiethammer11} and of this paper, and similarly of the convergence theorems of the two papers?

\medskip

There is a problem with the definition of the `integrated Rayleigh
principle' for a general function~$z$. Take the example of a Hilbert
space $H$ and a continuous semigroup generated by a non-negative
self-adjoint operator $A$, which solves the equation $\dot z = Az$ in
$H$. This is a gradient flow with $\langle z,Gy\rangle = (z,y)_H$ and
$\E(z) = \frac12 (z,Az)_H$. If $z(t)\not\in D(A^{1/2})$ at
$t>0$, then $\E(z(t))=\infty$ and $\mathrm D\E(z(t))$ is not well
defined. Another way of stating this is that the right-hand side
of~\eqref{def:RP} is not bounded from below, and its infimum equals
$-\infty$.  This in turn implies that even though~$\dot z$ might
minimize~\eqref{def:intRP} \emph{for fixed $z(\cdot)$}, in the
neighbourhood of that $z$ there exist perturbations $\widetilde z$,
arbitrary close to $z$, for which the infimum equals
$-\infty$. Therefore the formulation~\eqref{def:intRP} is very
unstable under perturbations of $z$.

\medskip

It is no coincidence that the expression~\eqref{def:intRP} is closely related to the functional $\A$. We can write
\begin{align}
\notag
\A(z) &= \E(z(T))-\E(z(0)) + \int_0^T \Biggl( 
\frac12 \langle \dot z, G(z)\dot z\rangle + \frac12 \Bigl\langle -\mathrm D\E(z), G(z)^{-1} (-\mathrm D\E(z))\Bigr\rangle\Biggr)\, \d t\\
&= \int_0^T \Biggl(\frac12 \langle \dot z, G(z)\dot z\rangle + \bigl\langle\mathrm D\E(z),\dot z\bigr\rangle \Biggr)\, \d t
+ \frac12 \int_0^T \Bigl\langle -\mathrm D\E(z), G(z)^{-1} (-\mathrm D\E(z))\Bigr\rangle\,\d t.
\label{eq:alt-form-A}
\end{align}
This form shows that Definition~\ref{def:GF} with the structure~\eqref{eq:29} is different from the integrated Rayleigh principle above in two ways:
%Comparing with the form of $J$ in~\eqref{def:J} yields some insight in the problem. The formulation in~\eqref{def:J} is different in two ways: 
first, the integral $\int_0^T \bigl\langle \mathrm D\E(z(t)),\dot z(t)\bigr\rangle\, \d t$ has been converted into the end point values $\E(z(T))-\E(z(0))$, and secondly, the addition of the dual dissipation potential $\psi^*(-\mathrm D\E(z))$ (the second term in~\eqref{eq:alt-form-A}) penalizes `non-regular' values of $z(t)$. 

Both of these changes appear to improve the robustness of the formulation. The addition of the dual potential has the effect of penalizing `unfavorable' choices for the function $z$; and the conversion of the cross term into end point values mitigates the effect of fast oscillations. 
The compactness results of Theorems~\ref{th:compactness-bounded} and~\ref{th:compactness} certainly suggest that Definition~\ref{def:GF} and~\eqref{eq:29} provide a useful basis for the analysis of these more general gradient flows. 

\bigskip

\textbf{Acknowledgements.}
A. Mielke was partially supported by the European Research Council
via "ERC-2010-AdG 267802" (AnaMultiScale).
The research of M. A. Peletier has received funding from the Initial Training Network
``FIRST'' of the Seventh Framework Programme of the European Community
(grant agreement number 238702).
G. Savar\'e has been partly supported by a grant from
MIUR for the PRIN08-project \emph{Optimal transport theory, geometric and functional inequalities and applications}.

\bibliography{Kramers}
%\bibliography{ref}
\bibliographystyle{alpha}
\end{document}